%% file: clearingenergymarkets_oms.tex
\documentclass[final,printer,a4paper]{gOMS2e}
\usepackage{fixltx2e}[2006/03/24]
\usepackage[T1]{fontenc}
\usepackage[utf8]{inputenc}
\usepackage{txfonts}  % CHOOSE Monotype Times ALTERNATIVE for OMS
\usepackage{amsmath}  % WARNING: amsthm is not compatible with OMS
\usepackage{eurosym}
\usepackage{algorithmic}
\usepackage{jmueller} % Modified version for OMS
\usepackage{graphicx}
\usepackage{threeparttable}
\usepackage[english]{babel}
\usepackage{url}
\usepackage{ifpdf}
\usepackage{float}%should be before hyperref
\ifpdf 
  \usepackage[pdfusetitle,colorlinks=true]{hyperref}
  \usepackage[figure]{hypcap}[2008/09/08]% hyperref to figures goes to the top of the image
\fi

\input{texDefinitions}

\input{commands}

\DeclareMathAlphabet{\mycal}{OMS}{mdbch}{m}{n} % DO NOT REMOVE THIS DEFINITION
\algsetup{indent=2em}%2em = qquad

\floatstyle{ruled}
\newfloat{algorithm}{thp}{lop}
\floatname{algorithm}{Algorithm}
\begin{document}
\doi{10.1080/10556788.2013.823544}
\issn{1029-4937}
\issnp{1055-6788}
\jvol{29} \jnum{1} \jyear{2014} \jmonth{}

\markboth{A. Martin, J.C. Müller, S. Pokutta}{Optimization Methods and Software}

\articletype{}%RESEARCH ARTICLE
\title{Strict linear prices in non-convex European day-ahead\\electricity markets}

\author{Alexander Martin$^{\rm a}$,
Johannes C. Müller$^{\rm a}$%
$^{\ast}$\thanks{$^\ast$Corresponding author. Email: Johannes.Mueller@fau.de},
and Sebastian Pokutta$^{\rm b}$\\\vspace{6pt}
$^{\rm a}${\em{Department of Mathematics,
FAU Erlangen-N\"urnberg,
Cauerstr. 11, D-91058 Erlangen, Germany}};
$^{\rm b}${\em{ISyE, Georgia Institute of Technology,
Groseclose 0205, 765 Ferst Dr, Atlanta, GA 30332, USA}};
%\\%
%{\em{Emails:
%Alexander.Martin@fau.de,
%Johannes.Mueller@fau.de,
%Sebastian.Pokutta@isye.gatech.edu
%}}
\\\vspace{12pt}\received{Received 9 February 2011; final version received 5 July 2013}}

\ifdefined\hypersetup
  \hypersetup{
    pdfauthor={Alexander Martin, Johannes C. Müller, Sebastian Pokutta},
    pdfkeywords={discrete optimization; MPEC; combinatorial auctions;
                 strict linear prices; electricity markets},
    linkcolor=black, citecolor=black, filecolor=black, urlcolor=black}
\fi

\maketitle

\begin{abstract}
The European power grid can be divided into several market areas where
the price of electricity is determined in a day-ahead auction. Market
participants can provide continuous hourly bid curves and
combinatorial bids with associated quantities given the prices. The
goal of our auction is to maximize the economic surplus of all participants
subject to quantity constraints and price constraints.
The price constraints ensure that no one incurs a loss. Only traders who
submitted a combinatorial bid might miss a not-realized profit.
The resulting problem is a large scale mathematical program with equilibrium
constraints (MPEC) and binary variables that cannot be solved
efficiently by standard solvers. We present an exact algorithm and
a fast heuristic for this type of problem. Both algorithms decompose
the MPEC into a master problem (a MIQP) and pricing subproblems (LPs).
The modeling technique and the algorithms are applicable to a wide
variety of combinatorial auctions that are based on mixed integer programs.
%157 words < 200 words
\bigskip
\begin{classcode}
90C11; %   	Mixed integer programming
90C33; %   	Complementarity and equilibrium problems and variational inequalities (finite dimensions)
%90C06; %   	Large-scale problems
91A46; %   	Combinatorial games
91B15; %   	Welfare economics
91B26. %   	Market models (auctions, bargaining, bidding, selling, etc.)
\end{classcode}
\bigskip
\end{abstract}

\section{Introduction}

In this article we present the market design that is currently in use by
most of the European electricity exchanges. The underlying
optimization problem is a welfare
maximizing combinatorial auction subject to price constraints.
For each tradable commodity a single price is determined to avoid
the same commodity of having different prices to different parties.
Even though the main
focus of the paper will be on solving the underlying problem, we also
explain why the market model evolved to its current state. Furthermore
we will discuss the differences to unit commitment models used in the US.

\subsection{Our contribution}
We provide a precise formulation of a real world optimization problem that
needs to be solved anew each day:
the market coupling problem between European day-ahead electricity exchanges.
The market is called a non-convex market because the feasible region of the
underlying optimization problem is non-convex. The non-convexity of the
feasible region is caused by the presence of binary decision variables.
In such a market it is difficult to find for each tradable commodity
a single price while ensuring that no participant incurs a loss.
The existence of such prices must be enforced by complementarity
conditions yielding a \emph{Mathematical Program with Equilibrium
Constraints (MPEC)} and binary variables. We present a deterministic algorithm yielding a provably
optimal solution and a fast
heuristic that exploits the model structure. Instead of solving the
large scale MPEC directly we decompose it into a Mixed Integer Quadratic
Program (MIQP) and a linear pricing problem. Empirical tests suggest that the
solutions determined by the heuristic are optimal in many cases.
We also ensure unique prices without changing
the economic surplus. Additionally, we introduce rules that allow
market participants to check whether the determined
prices satisfy necessary optimality conditions of the individual
participant's optimization problems.

\subsection{Outline} The structure of the article is as follows. In
Section~\ref{sec:preliminaries} we give a brief introduction to the European
electricity market and its features. In Section~\ref{sec:model} we
introduce the combinatorial optimization problem that is used to
determine welfare maximizing solutions subject to price constraints.
Analyzing this model we
obtain optimality conditions which are formulated in
Section~\ref{sec:optim-cond}. We then establish in
Section~\ref{sec:uniqueness-solution} the uniqueness of prices and
then we formulate a clearing heuristic and an optimal algorithm in
Section~\ref{sec:algorithm}. We
conclude with computational results in Section~\ref{sec:results}
and some final remarks in Section~\ref{sec:summary}.

\section{Background information}
In this section, we provide background information on linear prices and
on the institutional differences between US electricity markets and 
European electricity markets.

\subsection{Linear pricing schedules}
The European day-ahead electricity auctions are using linear pricing schedules.
We will briefly define the concept of pricing schedules, here in
the multidimensional setting that is a straightforward extension of 
\cite[p. 136]{tirole1988}.

\begin{definition}
We consider $m$ distinct commodities. A \emph{pricing schedule}
$T:\R^m\rightarrow\R$ is a map that returns the total amount of money 
to be paid by a consumer depending on his consumption vector
$q\in\R^m$. The schedule is called a \emph{linear pricing schedule} if
the map is linear, i.e., $T(q)=\p\tp q$. In this case $\p$ is called a 
\emph{linear price vector} and $\p_i$ is the price per unit for commodity $i$.
The definition is also applicable to producers if we use $-q_i$ 
consumption units to model $q_i$ production units of commodity $i$.
\end{definition}

\begin{definition}
The \emph{clearing condition} of a commodity $c$ is an equation
that ensures that the number of bought units of commodity $c$
is equal to the number of sold units of commodity $c$.
\end{definition}

\begin{definition}
A pricing schedule is a \emph{strict linear pricing schedule} if it is linear
and the number of commodities $m$ is equal to the number of clearing conditions
in the auction model \cite{vanvyve2011}.
In the electricity market, the clearing conditions are the flow conservation 
equations for each network node and time slot.
\end{definition}

In linear programming the strong duality theorem provides that
given a finite optimal solution to a primal linear program
there exists a finite optimal solution to the dual program.
The variables of the dual program are called dual variables and
can be interpreted as prices \cite[Chapter 12]{dantzig1963}.
In convex optimization strong duality holds if a constraint qualification
(e.g., Slater's condition) is satisfied. In this case, the dual variables 
can also be interpreted as prices. Sometimes the dual variables are
also called \emph{shadow prices} to the primal program \cite[Chapter 5]{boyd04}.

Now assume that $x$ is an optimal finite solution to a convex auction model
that maximizes the economic surplus. Let Slater's condition hold and let
$\p$ be optimal dual variables of the clearing conditions of all commodities.
Then $T(q)=\p\tp q$ defines a strict linear pricing schedule and $\p$ defines
a linear price vector. All participants are perfectly happy with the solution
$x$ and the prices $\p$, as their individual optimization problems are maximized
at these prices. Such a solution is called a competitive or Walrasian
equilibrium (cf. \cite{blumrosenNisan2007}).

\subsection{The US market}

We begin with a short introduction to electricity markets in the US.
In the US there exist several 
pool-based electricity markets, operated by \emph{independent
system operators (ISO)}. The ISOs organize day-ahead auctions where
both suppliers and consumers can submit bids. Using these bids the
auction determines prices for each network node and each hour of the
next day. Finally a unit commitment model computes a dispatch
that minimizes total production cost while ensuring system reliability.
%(cf. \href{http://www.pjm.com/sitecore\%20modules/web/~/media/training/core-curriculum/ip-gen-101/gen-101-two-settlement.ashx}{pjm.com}).
A standard reference for unit commitment models
is \cite{hobbs2001}. These models allow the participants to submit the cost structure and
production capacities of each power plant. In particular it is
possible to submit startup costs and minimal production capacities which are
modeled by using binary variables. This detailed information allows the ISOs
to ensure system reliability while minimizing the total production cost.
In many unit commitment models the demand is assumed to be a price inelastic bid
with a fixed given quantity.

Let us hypothetically assume a unit commitment problem can be described by a linear program (LP).
Then the optimal dispatch
would be determined by solving the LP and a strict linear pricing schedule is
given by the dual variables of the flow conservation equations of the
network nodes. As soon as binary variables are involved
we have to solve a mixed integer linear program (MIP) and the determination
of a reasonable strict linear pricing schedule is not as easy as in the LP case. 
Given an optimal solution to a MIP,
in general it is impossible to find a strict linear pricing schedule
where all participants neither incur a loss nor miss a not-realized profit.
Definition \ref{def non-realized profit} formally introduces the concept
of a not-realized profit. O'Neill et al.~\cite{oneill2005}
present a method that defines a linear pricing schedule for an
optimal solution to a MIP based unit commitment model.
The main idea is to treat binary actions like separate commodities.
Then we can find linear prices for all commodities, including the binary
actions. No participant incurs a loss and no participant misses a not-realized profit.
Each binary action can be attributed to a single market participant
and only this participant can carry out this binary action.
For this reason the price for the individual binary action can be
understood as an
individual compensation payment that is paid or received by each participant.
According to \cite{oneill2005}, the electricity auctions of
\emph{New York Independent System Operator} (\href{http://www.nyiso.com}{NYISO}) and
\emph{Pennsylvania-New Jersey-Maryland Interconnection} (\href{http://www.pjm.com}{PJM})
are using similar pricing methods.
Another technique that defines prices for a MIP solution is called convex hull pricing
(\cite{gribik2007}, \cite{wang2011}). This pricing schedule reduces the impact
of compensation payments but still requires the use of individual compensation payments.

In Europe the common practice is to avoid compensation payments and to
implement strict linear prices. This major difference renders methods
that are based on compensation payments inapplicable.

\subsection{The European market}
Similar to the US there exist several national and regional \emph{transmission system
operators (TSO)}\footnote{Nordic TSOs (NPS region):
  \emph{Elering, Energinet.dk, Fingrid, Litgrid, Statnett SF, Svenska Kraftn\"at}.
  Central Western Europe TSOs (CWE region): \emph{Amprion, Creos Luxembourg,
  Elia System Operator, Rte, Tennet TSO, TransnetBW}.
%Elering (Estonia), Statnett SF (Norway), Svenska Kraftn\"at (Sweden)
%Amprion (Germany), Elia System Operator (Belgium and Germany (50Hertz)), Rte (France)
%Tennet TSO (Netherlands and Germany), TransnetBW (Germany)
  (cf. \href{https://www.entsoe.eu/the-association/members/}{www.entsoe.eu})} in Europe.
The main difference to the ISOs in the US is that the European TSOs only
maintain system feasibility and reliability. They do not control the spot market
for electricity, but they determine network boundary conditions within which
they can guarantee the system feasibility and reliability. These boundary conditions
(e.g., available transmission capacities) are submitted to independent \emph{power
exchanges (PX)}\footnote{Nordic PX (NPS region): \emph{Nordpool Spot AS}.
Central Western Europe PXs (CWE region): \emph{APX Endex, Belpex, EPEX Spot}.}
who will run the financial spot market auctions subject to the network
constraints. Some power exchanges also split their market into several \emph{bidding areas}.
These bidding areas are often given by national borders or network bottlenecks.
For example the power exchange \emph{\epex} splits its market into the bidding areas
\emph{France, Germany/Austria, and Switzerland}.

\begin{figure}[tpb]
	\centering
	\begin{minipage}[c]{9cm}
\includegraphics[width=1.00\textwidth]{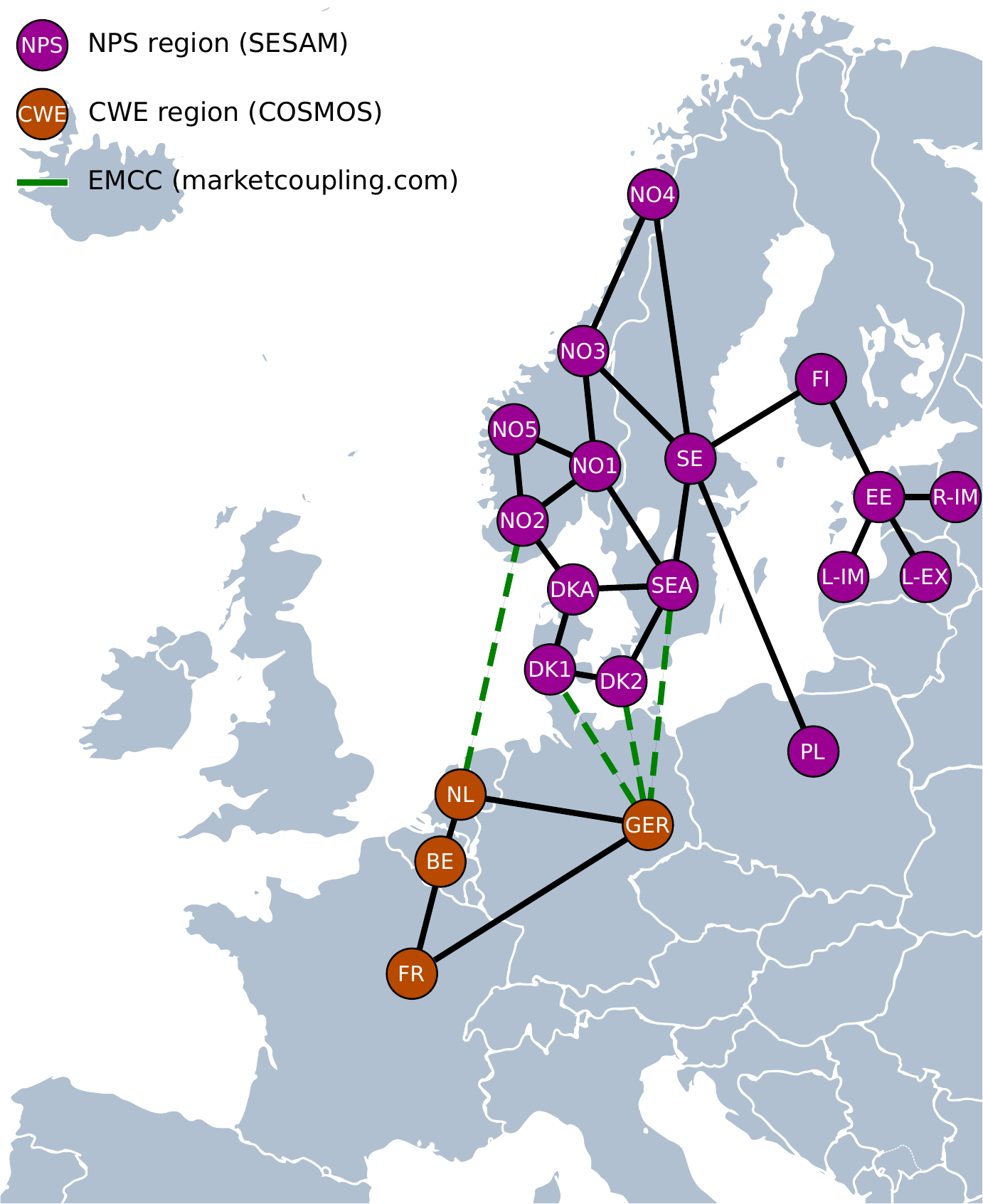}
	\centering
	\end{minipage}
	\caption{NPS and CWE region in 2011 (The Europe map is taken from Wikipedia)}
	\label{fig:network}
\end{figure}

On top of these power exchanges, \emph{market coupling systems}
collect all order book data from the power exchanges and network constraints from the TSOs and
compute welfare maximizing power flows between adjacent market areas and prices for the 
bidding areas of the power exchanges. These flows and price
signals are submitted back to the exchanges and TSOs. Finally the exchanges run their own
auction incorporating the cross-border flows and the price signals.
For example the \emph{SESAM} system couples the bidding areas of the
\emph{Nordic region (NPS region)} that consists of Denmark, Estonia, Finland, Lithuania, Norway,
and Sweden. The \emph{COSMOS} system couples the bidding areas of the 
\emph{Central Western European region (CWE region)} that consists of Belgium, France, Germany,
Luxembourg, and Netherlands. On top of these coupling systems, 
\emph{European Market Coupling Company (EMCC)} is computing the flows between the
NPS and the CWE region. SESAM is in use since 2007 and in Nov.\ 2009 EMCC started to
couple the NPS region with Germany. In Nov.\ 2010 COSMOS was launched and
from there on the EMCC system is coupling the NPS and the CWE region.\footnote{cf.
\href{http://www.nordpoolspot.com/About-us/History}{www.nordpoolspot.com/About-us/History},
\href{http://www.epexspot.com/en/market-coupling}{www.epexspot.com/en/market-coupling}, and
\href{http://www.marketcoupling.com/market-coupling/european-market}{www.marketcoupling.com}.}
Figure \ref{fig:network} shows the interconnected NPS and CWE region
in 2011. Note that these systems are distinct entities and guided by
different regulations while maintaining a certain set of compatibility
constraints as put forth in the European regulation of energy
markets.

All European power exchanges that we consider here are already existing and are using 
strict linear pricing schedules. Electricity at a specific bidding area at the specific
time it is being delivered is viewed as a commodity. Thus, for each bidding
area and each hour of the next day an electricity price (\eur/MWh) is computed, such
that no compensation payments are needed.
This will not be subject to change in the foreseeable future. 
The strict linear pricing schedule is determined by a welfare maximizing
optimization problem subject to quantity constraints and additional price constraints.
In order to perform a correct market coupling between such power exchanges a 
market coupling system needs to model the underlying power
exchanges exactly as they are. For this reason a European market coupling system
must also implement a strict linear pricing schedule and thus cannot implement
any other pricing schedule.

As of now, the TSOs in the NPS and the CWE region only submit upper and 
lower bounds for the flows and upper bounds for the change of flows.
Therefore a market coupling system that couples these two regions
only implements these basic transmission constraints.
However the algorithms that we present can be easily
adjusted to solve more complex transmission constrains, as long as
the constraints are linear or convex.
In the CWE region the introduction of so called \emph{power transmission 
distribution factor} matrices (PTDF) is planned for 2014 \cite{aguado2012flow}.

%%%%%%%%%%%%%%%%%%%%%%
%%%% Preliminaries
%%%%%%%%%%%%%%%%%%%%%%

\section{Bidding language in European day-ahead electricity auctions}
\label{sec:preliminaries}

Usually the auctioneer defines a \emph{bidding language} that allows the
participants for expressing their \emph{bids}, or in other words,
their preferences.
In a day-ahead electricity auction participants bid today for
electricity that they want to buy or sell on the following day.
They are allowed to submit their bids up to a predetermined point
in time. As soon as the submission deadline has passed 
the submitted bids cannot be changed anymore and the auctioneer
has to decide which bids to accept and which bids to reject.
Then the auctioneer publishes a pricing schedule that determines the
amount of money that a participant has to pay or receive if his bid
was accepted.

The day-ahead electricity auctions in the CWE and NPS region have a
common bidding language. The bidding language allows for submitting
four different types of bids: \emph{hourly bid curves},
\emph{block bids}, \emph{flexible bids}, and \emph{cross border trades}.
The number of different commodities the participants can bid for
is given by the number of market areas times the number of time slots
of the next day (e.g., 24 hours).
For this reason an entry of a strict linear price vector reflects
the price for one electrical power unit that is delivered at a specific
market area throughout a specific time slot.

In order to refer to a specific area or time slot we will use the following
notation: The set $A$ denotes the \emph{set of market areas} and $T$ denotes
the \emph{set of time slots}. An entry of a strict linear price vector $\p$
will be called \emph{price in area \areas\ at time \hours} and is denoted by
$\pat\in\R$.
We now define the four bid types in a formal way that suits the later exposition.

\subsection{Hourly bid curve}
A participant can submit hourly bid curves for every area and hour to
specify the optimal quantities that he is willing to buy or sell depending
on a given price. At first we will study a simple \emph{linear bid curve}
to understand the details of bid curves, then we will see that the concept of
linear bid curves can be used to express \emph{piecewise linear bid curves}.

\begin{figure}[tb]
	\centering
	\begin{minipage}[c]{\textwidth}
		\includegraphics[width=1.00\textwidth]{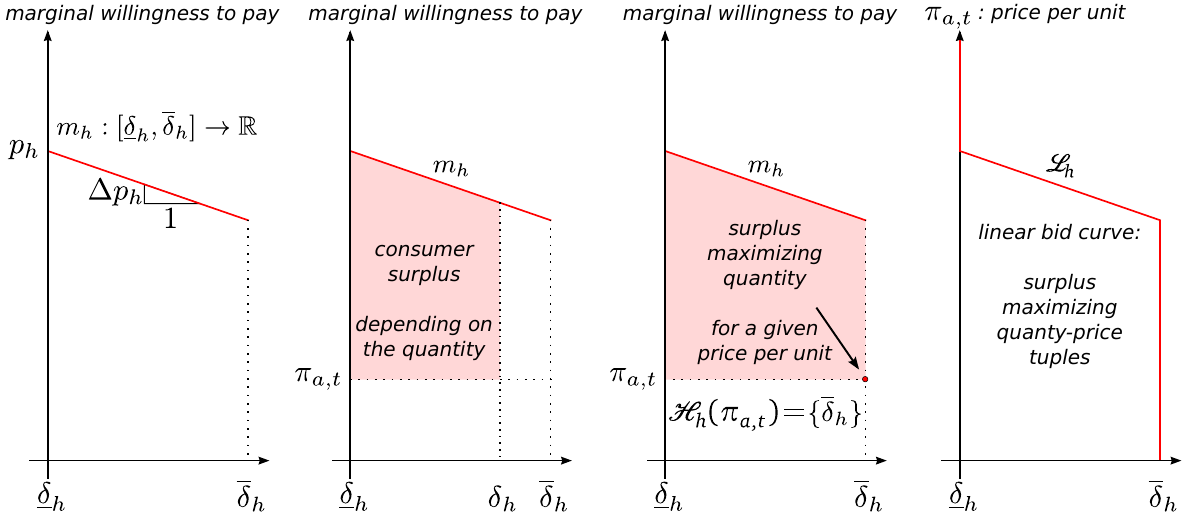}
		\centering
	\end{minipage}
	\caption{marginal willingness to pay, consumer surplus,
                 surplus maximizing quantities, and linear bid curve}
	\label{fig:willingness_to_pay}
\end{figure}
Assume that a participant submits a linear bid curve $h$
to buy electricity in area $a$ at time $t$.
He wants to buy at most $\dhub>0$ units and at least $\dhlb=0$ units.
Let his \emph{marginal willingness to pay} be given by an affine 
non-increasing function $m(u)=p_h-\dph u$.
Figure \ref{fig:willingness_to_pay} depicts his marginal willingness to pay.
The function valuates the price that the consumer is willing to pay
for each additional unit. If the function is constant ($\dph=0$), then
the price he is willing to pay for the first consumed unit is equal to
the price he is willing to pay for the last consumed unit.
His \emph{willingness to pay} %(also called \emph{valuation})
for the consumption of $\delta_h\in[\dhlb,\dhub]$
quantity units is given by $\int_0^{\delta_h}\phvonq(u)\,\dd u$.
If $m_h$ is constant, then this term reduces to $p_h\cdot\delta_h$. In this case
he would buy electricity if the price per unit is lower or equal to $p_h$.
We will now assume that the participant is interested in maximizing his
\emph{consumer surplus}.

\begin{definition}
The \emph{consumer surplus} of a participant is his willingness to pay for
the consumption of $q$ units minus the total amount of money to be paid for
$q$ units (cf. \cite[p. 8]{tirole1988}).%
%\footnote{\cite{tirole1988} Page 8: ``Dupuit (1844) contended that this area
%was a measure of what the consumers would be willing to pay in excess of
%what they already spend ($\pat\cdot\delta_h$) for the right to consume $\delta_h$
%units of the good.''}
\end{definition}

Let the price per unit in area $a$ at time $t$ be exogenously given by the
auctioneer and let it amount to $\pat$ currency units per quantity unit.
The participant has only one decision variable: the quantity $\delta_h$
that he wants to buy. He can maximize his consumer surplus by solving the 
following parameterized optimization problem
(\pat\ is an exogenously given parameter).
\begin{align*}
\max\quad&\int_0^{\delta_h}\phvonq(u)\,\dd u - \pat\delta_h
\tag{QP-Hourly}\label{QP-Hourly-1}\\
\st \quad&\dhlb\le\delta_h\leq\dhub
\end{align*}
\begin{definition}
Let $m_h:[\dhlb,\dhub]\rightarrow\pricerange$ be an affine non-increasing function
and let $0\in[\dhlb,\dhub]$. The set of \emph{surplus maximizing quantities}
depending on the price is given by
\[
\mycal{H}_h(\pat):=
\argmax_{\dhlb\le\delta_h\leq\dhub}
\int_0^{\delta_h}\phvonq(u)\,\dd u - \pat\delta_h.
\]
Note that $\mycal{H}_h(\pat)$ is equal to $[\dhlb,\dhub]$ if $m_h$ is
constant and $\pat=p_h$. The two-dimensional set
\[
\mycal{L}_h:=\left\{\vect{\delta_h,\pat}\mid\pat\in\R\und \delta_h\in\mycal{H}_h(\pat)\right\}
\]
is called a \emph{linear bid curve}. It is called linear, because it contains
exactly one linear segment that indicates a change in the quantity.
Figure \ref{fig:willingness_to_pay} illustrates both definitions.
\end{definition}
In Section \ref{sec:model} we will study the optimality conditions that characterize
the surplus maximizing quantities; for now, it suffices to be aware of the
structure of linear bid curves as depicted in Figure \ref{fig:willingness_to_pay}.

\begin{figure}[tb]
	\centering
	\begin{minipage}[c]{\textwidth}
		\includegraphics[width=1.00\textwidth]{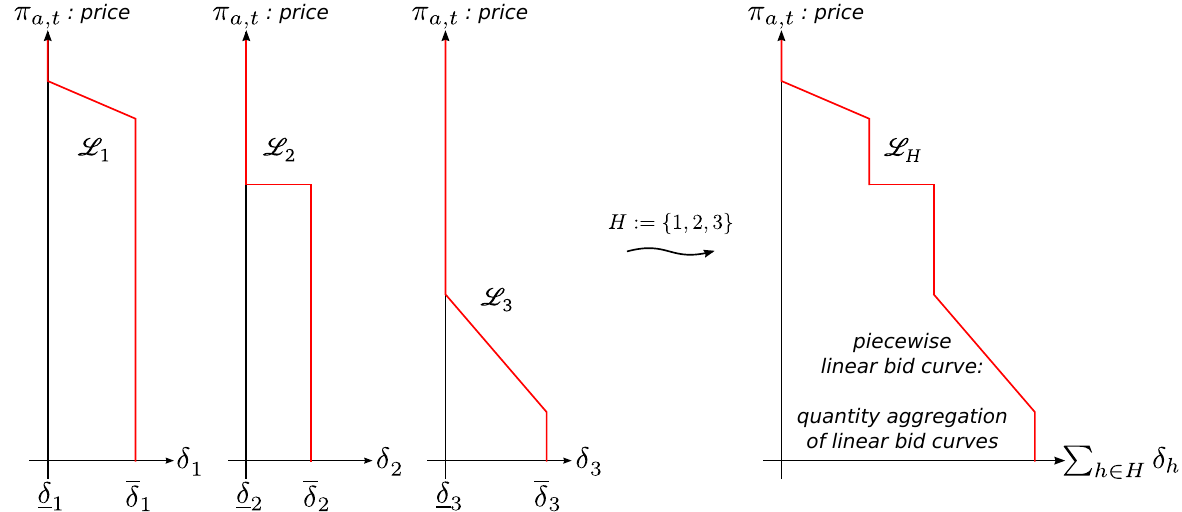}
		\centering
	\end{minipage}
	\caption{construction of piecewise linear bid curves with the help of linear bid curves}
	\label{fig:bid_curve_aggregation}
\end{figure}
If a participant wants to submit a more complex bid curve he can simply
submit a set $H$ of linear bid curves to construct a piecewise linear bid curve
$\mycal{L}_H$. The following definition formally introduces piecewise linear
bid curves and Figure \ref{fig:bid_curve_aggregation} provides an example.
\\
\begin{definition}
Let $H$ be a finite set of linear bid curves and
let $m_h:[\dhlb,\dhub]\rightarrow\pricerange$ be an affine
non-increasing function and $0\in[\dhlb,\dhub]$ for all $h\in H$.
The two dimensional set
\[
\mycal{L}_H:=\left\{\vect{\sum_{h\in H}\delta_h,\pat}
\mid\pat\in\R\und \delta_h\in\mycal{H}_h(\pat)\fuer h\in H\right\}
\]
is called a \emph{piecewise linear bid curve}.
Let $\Hat$ be the set of all linear bid curves in area $a$ at time $t$.
The piecewise linear bid curve $\mycal{L}_{\Hat}$ is called
the \emph{aggregated bid curve} of area $a$ at time $t$.
\end{definition}

\begin{figure}[tb]
	\centering
	\begin{minipage}[c]{\textwidth}
		\includegraphics[width=1.00\textwidth]{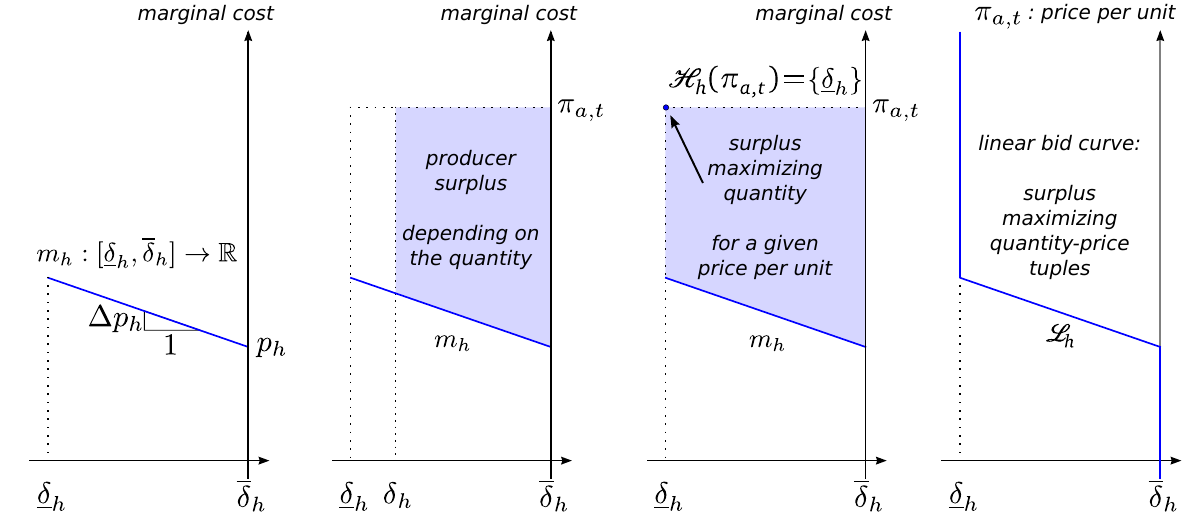}
		\centering
	\end{minipage}
	\caption{marginal cost, producer surplus, surplus maximizing quantity,
                 and linear bid curve}
	\label{fig:marginal_cost}
\end{figure}

Now suppose that a participant wants to sell electricity in area $a$ at time $t$.
He is able to generate at most
$|\dhlb|$ units and at least $\dhub=0$ units. The maximal supply quantity is modeled
by a negative number $\dhlb<0$. This allows us for using the definitions from above,
such that we do not need to distinguish between demand and supply.
A positive quantity indicates demand and a negative quantity indicates supply.
Like above, we assume that the \emph{marginal cost curve}
$m_h:[\dhlb,\dhub]\rightarrow\pricerange$ is an affine non-increasing function: 
the last supplied unit (unit $\dhlb$) is the most expensive one, the first supplied unit is
the cheapest one. The \emph{cost} for generating $\delta_h\in[\dhlb,\dhub]$ units
is given by $\int_{\delta_h}^0\phvonq(u)\,\dd u$.
% +++ das ist richtig, auch wenn man aufgrund der economies of scale
% etwas anderes erwartet hätte.
% die ersten nachgefragen megawatts wird man am liebsten mit der
% billigsten produktionseinheit herstellen: atomkraft.
% wenn das atomkraftwerk am limit arbeitet, dann müssen
% z.B. gasturbinen hinzugeschaltet werden. deren marginale kosten
% sind höher als die von atomkraftwerken.
\begin{definition}
The \emph{producer surplus} of a participant is the amount of money to
be received for producing $q$ units minus the cost for generating $q$ units.
Figure \ref{fig:marginal_cost} provides an example.
\end{definition}
In Equation \eqref{eq producer surplus} it can be shown that the producer and consumer
surplus can be expressed by the same formula. For this reason the
producer can maximize his producer surplus by solving \eqref{QP-Hourly-1}.
\begin{equation}\label{eq producer surplus}
\pat(-\delta_h) - \int_{\delta_h}^0\phvonq(u)\,\dd u= 
\int_0^{\delta_h}\phvonq(u)\,\dd u - \pat\delta_h
\end{equation}

\subsubsection{Notation used for hourly bid curves}
As we have seen in the previous section, we do not need to model piecewise linear
bid curves directly. It is sufficient to treat each underlying linear bid curve
separately. For this reason we only need one set of linear bid curves for each
area and time to model all piecewise linear bid curves that were submitted
by the participants.
\begin{align*}
  &H&&\text{set of all linear bid curves}\\
  &H_{a,t}\subseteq H&&\text{set of all linear bid curves in area $a$ at time $t$}\\
  &\dhlb\le0\le\dhub&&\text{quantity bounds of the linear bid curve $h\in H$}\\
  &m_h:[\dhlb,\dhub]\rightarrow\pricerange&&
\text{affine non-increasing marginal willingness to pay / marginal cost for $h\in H$}\\
  &\delta_h\in[\dhlb,\dhub]&&
\text{decision variable: bought quantity if positive / sold quantity if negative}
\end{align*}

\subsection{Block bids}
Block bids enable traders to buy or sell electricity in an area for several (not
necessarily consecutive) time slots with a single bid.
A block bid can either be executed entirely or it is not executed at all.
This condition is called \emph{fill-or-kill} condition.
%\footnote{\url{http://www.sec.gov/answers/fokord.htm}}
Suppose that a participant wants to buy electricity in area $a$ and wants to use
a block bid $b$. Then for each hour \hours\ he has to specify the quantity
$q\bt$ that he wants to buy or sell. Positive quantities indicate that
he wants to buy electricity and negative quantities indicate that he wants
to sell electricity. Currently it is only allowed to transmit either a pure
buy bid with $q_{b,t}\ge0$ or a pure sell bid with $q_{b,t}\le0$ for all \hours{}
(cf. \cite{eex2013rules}).
Let $\beta_b\in\bin$ be the binary variable that models
whether the bid \(b\) is $(\beta_b=1)$ or is not $(\beta_b=0)$ executed.
Then $\sum_\hours q\bt\beta_b$ is the total bought (or sold) quantity.
Let $p_b$ be his constant \emph{marginal willingness to pay},
%\footnote{At stock and power exchanges this price is called \emph{limit price}%
%, cf.\ \url{http://www.sec.gov/answers/limit.htm}
%}.
%\footnote{In \cite[p. 318]{mascolell1995} the marginal willingness to pay
%    is called \emph{marginal benefit} $\Phi'_b(x)=p_b$ and the term 
%    $\Phi_b(x)=\int_0^x\Phi'_b(u)\dd u$ is called \emph{benefit}. 
%    Note that $\Phi_b(\sum_\hours q\bt\beta_b)=p_b\cdot\sum_\hours q\bt\beta_b$}
i.e., for each bought electricity unit he is willing to pay $p_b$
currency units. His \emph{willingness to pay} for the consumption of
$\sum_\hours q\bt\beta_b$ units is given by $p_b\cdot\sum_\hours q\bt\beta_b$.
Let the price $\pat$ in area $a$ and time $t$ be given, then
the price-parameterized surplus maximization problem for bid $b$ 
looks as follows:
\begin{align*}
\max\quad&p_b\sum_\hours q\bt\beta_b - \sum_\hours \pat q\bt\beta_b
\label{MIP-BlockBid}\tag{MIP-BlockBid}
\\\st\quad&\beta_b\in\bin
\end{align*}
The auctioneer cannot guarantee that at the end of the auction the 
binary variable $\beta_b$ is a surplus maximizing solution to
\eqref{MIP-BlockBid}, because this is simply impossible in general.
\begin{example}\label{ex paradoxically rejected}
Let the set of areas and the set of hours be singletons: $A=\{a\}$, $T=\{t\}$.
Bid $b$ is a sell block bid with quantity $q_{b,t}=-1$ MWh and
marginal cost $p_b=1$ \eur/MWh. Bid $c$ is a buy block bid with quantity 
$q_{c,t}=2$ MWh and marginal willingness to pay $p_c=2$ \eur/MWh. The 
two orders cannot be matched, because of the different quantities. The only
feasible solution is to reject both bids, that is, $\beta_b=\beta_c=0$.
If we publish a price that is strictly smaller than $2$ \eur/MWh, then bid 
$c$ is wondering why it is rejected. If we publish a price that is 
strictly greater than $1$ \eur/MWh, then bid $b$ is wondering why it is rejected.
Therefore, there exists no price where everyone is perfectly happy:
there is no price $\pat$ such that $\beta_b$ maximizes \eqref{MIP-BlockBid}$_b$
and $\beta_c$ maximizes \eqref{MIP-BlockBid}$_c$.
\end{example}

Instead of ensuring the individual surplus maximization of a block bid,
the auctioneer ensures that the surplus is non-negative
if the block bid is executed. In other words, a buy block bid can
only be executed if the willingness to pay is greater or equal to
the amount of money to be paid:
\begin{equation}\label{eq block price condition}
p_b\sum_\hours q\bt\beta_b - \sum_\hours \pat q\bt\beta_b \ge 0.
\end{equation}
If the block bid is a sell block bid, then all quantities are smaller or equal 
to zero: $q\bt\le0$. In this case $p_b$ is the \emph{marginal cost} of production and
equation \eqref{eq block price condition} ensures that the sell block bid
can only be executed if the money to be received for generating electricity
is greater or equal to the cost for generating the electricity:
$- \sum_\hours \pat q\bt\beta_b \ge -p_b\sum_\hours q\bt\beta_b.$
\begin{definition}\label{def non-realized profit}
A block bid $b$ \emph{incurs a loss} if equation
\eqref{eq block price condition} is violated. A block bid $b$
\emph{misses a not-realized profit} if it is rejected even though
\begin{equation}\label{eq block price profit}
p_b\sum_\hours q\bt - \sum_\hours \pat q\bt > 0,
\end{equation}
holds. Block bids that miss a not-realized profit, are called
\emph{paradoxically rejected bids (PRB)}, because they are rejected
even though from the local point of view of the owner of the block the rejection
is not a surplus maximizing solution to \eqref{MIP-BlockBid}.
\end{definition}

In Example \ref{ex paradoxically rejected} either block $b$ or block $c$ is
paradoxically rejected, showing that paradoxically rejected blocks are inevitable.

The trader further has the possibility to \emph{link} certain block
bids so that the execution of some block $b_1 \in B$ implies the execution
of some other block $b_2 \in B$.

\subsubsection{Notation used for block bids}
\ \vspace*{-1.7\baselineskip}% The space makes this code latex version independent
\begin{align*}
  &B&&\text{set of all block bids}\\
  &B_a\subseteq B&&\text{set of all block bids in area $a$}\\
  &L\subseteq B\times B&&\text{%the set of all block links:
if $(b,c)\in L$, then $b$ can only be executed if also $c$ is executed}\\
  &q_{b,t}\in\Q&&
\text{quantity to be bought ($\ge0$) / sold ($\le0$) in hour \hours\ by block bid $b$}\\
  &p_b\in\pricerange&&\text{constant marginal willingness to pay / marginal cost of block bid $b$}\\
  &\beta_b\in\bin&&\text{decision variable: rejected ($0$) / executed ($1$)}
\end{align*}

\subsection{Flexible bids}
With flexible bids a participant can buy or sell electricity in a specific area in 
exactly one hour without specifying the hour in which the bid should be executed. 
Flexible bids must satisfy the \emph{fill-or-kill} condition as well,
i.e., a partial execution is not feasible. Suppose that a participant
wants to buy $q_f>0$ units in area $a$ in a not specified time slot. Let his
\emph{marginal willingness to pay} amount to $p_f$ currency units per quantity unit. 
We need to introduce one binary variable $\pft\in\bin$ for each time slot \hours.
If $\pft=1$, then the flexible bid $f$ is executed in time slot $t$. The bid
may only be executed in at most one time slot: $\sum_\hours \pft\le1$. The
\emph{willingness to pay} for the total consumed quantity of $q_f\sum_\hours \pft$
units is given by $p_f q_f\sum_\hours \pft$. Now we can write down
the price-parameterized surplus maximization problem of a flexible bid.
\begin{align*}
\max\quad&p_f q_f\sum_\hours \pft - \sum_\hours\pat q_f\pft
\label{MIP-FlexBid}\tag{MIP-FlexBid}\\
\st\quad&\sum_\hours \pft\le1\\
        &\pft\in\bin\qquad\fuer\hours
\end{align*}
In general it is impossible to ensure that at the end of the auction all decision
variables $\pft$ will be surplus maximizing for \eqref{MIP-FlexBid}, but it
is possible to ensures that the surplus is
non-negative if the flexible bid is executed:
\begin{equation}\label{eq flex price condition}
  p_f q_f\sum_\hours \pft - \sum_\hours\pat q_f\pft \ge 0.
\end{equation}
If this inequality is violated, then the flexible bid \emph{incurs a loss}.
At the end of the auction some flexible bids might \emph{miss a not-realized profit},
i.e., from the local point of view of the owner of a flexible bid the current solution
is not necessarily surplus maximizing for \eqref{MIP-FlexBid}.

\subsubsection{Notation used for flexible bids}
\ \vspace*{-1.7\baselineskip}% The space makes this code latex version independent
\begin{align*}
  &F&&\text{set of all flexible bids}\\
  &F_a\subseteq F&&\text{set of all flexible bids in area $a$}\\
  &q_f\in\Q&&\text{quantity to be bought ($>0$) / sold ($<0$) by flexible bid $f$}\\
  &p_f\in\pricerange&&\text{constant marginal willingness to pay / marginal cost of flexible bid $f$}\\
  \pft&\in\{0,1\}&&\text{decision variable: }\pft=\begin{cases}
		1&\text{flexible bid \flex\ is executed in hour \hours}\\
		0&\text{no execution}
\end{cases}
\end{align*}

\subsection{Cross-border trades}
Trading electricity across areas is called a \emph{cross-border trade}.
The regulatory body stipulates that in European day-ahead electricity auctions, only transmission
system operators are allowed to perform cross border trades. The 
cross-border trader buys electricity in low price areas and sells electricity
in high price areas. Cross-border trades are only possible if there are 
\emph{interconnectors} connecting the areas. The latter will be denoted 
by $C\subseteq A\times A$, where $A$ is the set of areas. An interconnector
$c=\rs\in C$ is modeled by a directed arc, whereas $r$ denotes the source
and $s$ denotes the sink. The transmitted electricity on interconnector
\cons\ in time slot \hours\ is denoted by a signed variable $\tauct$.
We use the convention that a positive flow indicates a transfer from 
the source $r$ to the sink $s$, and a negative flow indicates a transfer
from $s$ to $r$. There are upper and lower bounds $[\taulb,\tauub]$ on 
the transmission quantity. These bounds are called \emph{available
transmission capacity} (ATC). Also the change of flow between two consecutive
hours can be limited by an upper bound $\tauramp$, called \emph{ramp rate}.
Now suppose that there is a positive flow $\tauct>0$ on interconnector
$c=\rs$ in hour $t$, then electricity is bought in area $r$ and 
sold in area $s$. The \emph{marginal cost} for buying electricity is
$\prt$ and the \emph{price per unit} to be received for selling electricity
is $\pst$. The surplus of a cross-border trade is called \emph{congestion rent}
and can be maximized by solving the following price-parameterized
optimization problem.
\begin{align*}
\max\qquad&
	\sum_\hours(\pst-\prt)\tauct\tag{LP-TSO}\label{LP-TSO-1}
\\
\st\qquad
	&\forall\hours\qquad
\taulb\le\tauct\le\tauub\\
	&\forall\hours\qquad
-\tauramp\le\tauct-\tau_{c,t-1}\leq\tauramp
\end{align*}
The parameter $\tau_{c,-1}$ is the given flow of the last hour of the previous 
day.

At the end of the auction the flows will be congestion rent maximizing,
but we need to keep in mind that the cross-border trader is
considered a \emph{price taker}, i.e., the price is exogenously given by the auctioneer.
We will later see that in this case the optimality conditions of cross-border
trades imply that prices only diverge, if a transmission constraint is active.

\subsubsection{Notation used for cross-border trades}
\ \vspace*{-1.7\baselineskip}% The space makes this code latex version independent
\begin{align*}
  C&\subset A\times A&&\text{set of all interconnectors}\\
  C_a^+&=\{\rs\in C\,|\,r = a\}&&\text{set of all interconnectors starting in $a$}\\
  C_a^-&=\{\rs\in C\,|\,s = a\}&&\text{set of all interconnectors ending in $a$}\\
  \tau_{c,-1}&\in\Q&&\text{flow on $c$ in the last time slot of the
    previous day}\\
 \tauramp&\in\Q^+_0\cup\{\infty\}&&
\text{ramp rate; we put $\tauramp=\infty$ if we have no ramping on $c$}\\
  \tauct&\in[\taulb,\tauub]&&\text{decision variable: flow on \cons\ in time \hours\ within bounds (ATC)}\\
\end{align*}

%%%%%%%%%%%%%%%%%%%%%%
%%%% Model 
%%%%%%%%%%%%%%%%%%%%%%

\section{Optimization problem of the auctioneer}
\label{sec:model}
In this section we will present the complete MPEC with binary variables
that describes the auctioneers optimization problem. The aim of
the auctioneer is to maximize the total economic surplus subject to
the clearing condition and the constraints given by the submitted bids.
We first introduce the overall model structure and then discuss the
objective function as well as each constraint in detail.

\begin{figure*}
\begin{align*}
 \max\enskip
	&\text{\it (\surplus)}
	&&\sum_{h\in H}\int_0^{\delta_h}\phvonq(u) \dd u
		+\sum_{\substack{\blocks\\\hours}} p_b\,q_{b,t}\,\beta_b
		+\sum_{\substack{\flex\\\hours}} p_f\,q_f\,\pft
		\tag{MPEC}\label{QPFlow}\\
\st\enskip
&\text{\it (quantity constraints)}\\
%[\p]\enskip
&\forall a\in A,\,t\in T	&&
		\sum_\hourly\delta_h +\sum_\blocka q_{b,t}\beta_b
		 + \sum_\flexa q_f \pft
		=\sum_\cin \tau_{c,t} - \sum_\cout \tau_{c,t}
		\tag{\ref{nb:clearing}}\\
%[\rho]
\enskip&\forall c\in C,\,t\in T	&&
 		-\tauramp\leq\tau_{c,t} - \tau_{c,t-1}\leq\tauramp
 		\tag{\ref{nb:ramprate}}\\
&\forall (b,c)\in L&&	\beta_b\leq\beta_c
		\tag{\ref{nb:blocklinks}}\\
&\forall f\in F&&	\sum_\hours \pft\leq 1
		\tag{\ref{nb:flexphi}}\\
&\forall\blocks&&	\beta_b\in\{0,1\}\\
&\forall\flex,\hours&&	\pft\in\{0,1\}\\
%[\mu]\enskip
&\forall\cons,\hours&&	\tau_{c,t}\in[\taulb,\tauub]\\
%[v]\enskip
&\forall h\in H&&	\delta_h\in[\dhlb,\dhub]\\
&\text{\it (price constraints)}\\
&&&\text{\it strict linear prices:}\\
&\forall\areas,\hours&&	\pat\in\pricerange\\
&&&\text{\it non-negative surplus of block and flexible bids:}\\
&\forall\areas, \blocka&&
		\beta_b\sum_\hours(p_b - \pat)\,q_{b,t}\geq0
		\tag{\ref{nb:blockpreisbedingung}}\\
&\forall\areas, \flexa, \hours&&
		\pft(p_f-\pat)\,\sgn(q_f)\geq0
		\tag{\ref{nb:flexpreisbedingung}}\\
&&&\text{\it optimality conditions of cross-border flow, \eqref{LP-TSO-1}:}\\
&\forall c=\rs\in C,\hours&&
%		(\muub\ct -\mulb\ct) +(\rhoub\ct -\rholb\ct) +[(\rholb\ctp -\rhoub\ctp)]+\prt-\pst=0
		(\muub\ct -\mulb\ct) +(\rhoub\ct -\rholb\ct) +[(\rholb\ctp -\rhoub\ctp)]=\pst-\prt
		\tag{\ref{nb:fpc}}\\
&\forall\cons,\hours&&
		(\tauct-\tauub)\;\muub\ct=0
		\tag{\ref{nb:fpc}}\\
&\forall\cons,\hours&&
		(-\tauct+\taulb)\;\mulb\ct=0
		\tag{\ref{nb:fpc}}\\
&\forall\cons,\hours&&
		(\tauct-\tau_{c,t-1}-\tauramp)\;\rhoub\ct=0
		\tag{\ref{nb:fpc}}\\
&\forall\cons,\hours&&
		(\tau_{c,t-1}-\tauct-\tauramp)\;\rholb\ct=0
		\tag{\ref{nb:fpc}}\\
&\forall\cons,\hours&&
		\muub\ct,\;\mulb\ct,\;\rhoub\ct,\;\rholb\ct\geq 0
		\tag{\ref{nb:fpc}}\\
&&&\text{\it optimality conditions of linear bid curves, \eqref{QP-Hourly-1}:}\\
&\forall\areas,\hours,\hourly&&
		-\phvonq(\delta_h)+\pat+\vub_h-\vlb_h=0
		\tag{\ref{nb:fillingcondition}}\\
&\forall\areas,\hours,\hourly&&
		(\delta_h-\dhub)\;\vub_h = 0
		\tag{\ref{nb:fillingcondition}}\\
&\forall\areas,\hours,\hourly&&
		(-\delta_h+\dhlb)\;\vlb_h = 0
		\tag{\ref{nb:fillingcondition}}\\
&\forall\areas,\hours,\hourly&&
		\vub_h,\;\vlb_h\geq 0
		\tag{\ref{nb:fillingcondition}}
\end{align*}
\end{figure*}
The primal variables of the model \eqref{QPFlow} are $\beta,\phi,\tau,\delta$ and
$\p,\rholb,\rhoub,\mulb,\muub,\vlb,\vub$. The former are used to model the executed 
quantities and the latter are used to model the price constraints. In the next
section we show that the primal variables $\rholb,\rhoub,\mulb,\muub,$ and $\vlb,\vub$
in this model are used to model the Karush-Kuhn-Tucker optimality conditions (KKT)
of the price-parameterized surplus 
maximization problems of cross-border trades \eqref{LP-TSO-1} and linear 
bid curves \eqref{QP-Hourly-1}. A short introduction to the 
KKT optimality conditions can be found in \cite[Chapter 5.5.3]{boyd04}.
Because of the KKT conditions, some of the price constraints are non-convex.
The bar below or above the variables $\rholb,\rhoub,\mulb,\muub,\vlb,\vub$
indicates that the variable is associated with a lower or upper quantity bound of a 
quantity constraint of \eqref{LP-TSO-1} or \eqref{QP-Hourly-1}.
For example the variable $\vub_h$ is associated with the upper quantity bound 
$\delta_h\le\dhub$ of \eqref{QP-Hourly-1}. We will make this precise below.
As the auctioneer needs to compute surplus maximizing
quantities to the surplus maximization problems of linear bid curves and
cross-border trades, the dual variables of these subproblems
are primal variables in the auctioneers optimization problem.
We need to explicitly add the optimality conditions (price conditions) to the model,
because they are not necessarily automatically fulfilled
in an optimal solution.

An important property of the model is that the determined strict linear prices are 
consistent with the strict linear prices that we know from convex auctions: Assume 
that we solve a relaxed model that only contains the quantity constraints
(except for constraint \nbref{nb:blocklinks}) and assume that its QP-relaxation admits
an integral optimal solution, with respect to the binary variables.
Then we can find strict linear prices to this optimal
solution that satisfy all previously mentioned price conditions. In this case, all price 
conditions are dispensable and could be omitted.
This assertion can be verified by writing down the KKT conditions
of the QP-relaxation. Van Vyve proves this assertion in Proposition 3 in \cite{vanvyve2011}.
%if we choose the set of rejected bids $J_0$ to be the empty set.
The same result in a similar context is also provided in 
\cite[Corollary 11.16]{blumrosenNisan2007}.

\subsection{Constraints}\label{section:constraints}
\renewcommand{\theenumi}{\text{C\arabic{enumi}}}
\renewcommand{\labelenumi}{(C\arabic{enumi})}
We now introduce the constraints mentioned above and explain them in detail. 
\begin{enumerate}

\item{\textbf{Clearing condition.}}\label{nb:clearing} This
  constraint ensures the identity of executed net demand and net
  import. The constraint is similar to the classical
  \emph{conservation of flow} constraint typically present in network
  flow models. In every area and every hour the executed net demand must coincide with the net import:
\[
 \forall a\in A,\,t\in T	\qquad
		\sum_\hourly\delta_h +\sum_\blocka q_{b,t}\beta_b
		 + \sum_\flexa q_f \pft
		=\sum_\cin \tau_{c,t} - \sum_\cout \tau_{c,t}\;.
\]

\item{\textbf{Ramping condition.}}\label{nb:ramprate} This
  condition incorporates the ramping requirements on the
  interconnectors by limiting the change between two consecutive time
  slots accordingly. From one hour to the next hour the flow on
  interconnector $c$ may only change by $\tauramp$ units:
\[
 \forall c\in C,\,t\in T		\qquad
 |\tau_{c,t} - \tau_{c,t-1}|\leq\tauramp\;.
\]

\item{\textbf{Block bid links.}}\label{nb:blocklinks} To model the
  linking of block bids we add the following constraint. If $(b,c)\in L$, then block $b$ may only be executed if block $c$ is executed:
\[
 \forall (b,c)\in L\qquad	\beta_b\leq\beta_c\;.
\]

\item{\textbf{Flexible bid execution at most once.}}\label{nb:flexphi}
  A flexible bid $f$ can be executed in at most one time slot: 
\[
 \forall f\in F\qquad	\sum_\hours \pft\leq 1\;.
\]

\item{\textbf{Block price condition.}}\label{nb:blockpreisbedingung}
  In principle it is possible to obtain an
  equally high or even higher economic surplus by requiring altruistic behavior
  from some of the participants. However this is not in line with
  general economic principles as no participant is willing to sell at
  a loss. Therefore, a block bid can only be executed if it does not
  incur a loss:
\begin{align*}
 \forall\areas,\blocka\qquad	\beta_b\sum_\hours(p_b - \pat)\,q_{b,t}\geq0\;.
\end{align*}
  As mentioned before the converse cannot be guaranteed:
  some block bids might miss a not-realized profit at
  the end of the auction. In other words, they might be rejected,
  even though the prices would allow the execution.

  The block price condition is a quadratic non-convex constraint that is difficult to handle,
  but it is easy to linearize it if we introduce a sufficiently large constant $M_b$ for
  each \blocks. We know that $\beta_b$ is a binary variable and we assume
  that $\pat\in[\pmin,\pmax]$. Now we can rewrite the constraint as follows
  (cf. \cite[Chapter 26-3-I(g)]{dantzig1963}):
\[
  \sum_\hours(p_b - \pat)\,q_{b,t}\geq M_b(1-\beta_b)\quad\text{with}\quad
  M_b := \sum_\hours\min\{(p_b - \p)\,q_{b,t}\mid\p\in[\pmin,\pmax]\}\;.
\]

\item{\textbf{Flexible bid price condition.}}\label{nb:flexpreisbedingung}
A flexible bid is only executed if it does not incur a loss:
\begin{align*}
\forall\areas,\flexa,\hours\qquad	\pft(p_f-\pat)\,\sgn(q_f)\geq0,
\end{align*}
where \(\sgn(.)\) denotes the sign function. This formulation is equivalent
to equation \eqref{eq flex price condition}. Similar to block bids,
also flexible bids might miss a not-realized profit, i.e., they might
get rejected, even though the prices would allow the execution.

\item{\textbf{Flow price condition.}}\label{nb:fpc}
This condition represents the optimality condition of the cross-border
traders optimization problem. The regulators force them to be price takers
such that the prices are externally given by the auctioneer and the
only remaining decision variable is the flow $\tauct$. The individual
optimization problem is to maximize the congestion rent subject to the
available transmission capacity and the ramp rate:
\begin{align*}
\max\qquad
		&\sum_\cons
				\sum_\hours(\pst-\prt)\tauct
				\tag{LP-TSO}\label{LP-TSO}\\
\st\qquad
		&[\muub]\qquad
\forall c\in C,t\in T
\qquad \tauct\leq \tauub&&\text{ATC}\\
		&[\mulb]\qquad 
\forall c\in C,t\in T
\qquad -\tauct\leq-\taulb&&\text{ATC}\\
		&[\rhoub]\qquad
\forall c\in C,t\in T
\qquad \tauct-\tau_{c,t-1}\leq\tauramp&&\text{ramp rate}\\
		&[\rholb]\qquad
\forall c\in C,t\in T
\qquad \tau_{c,t-1}-\tauct\leq\tauramp&&\text{ramp rate}
\end{align*}
The variables in squared brackets are the dual variables to the primal
constraints. By duality, essentially an optimal flow implies that a price difference
between two adjacent areas implies that at least one transmission constraint is
active. %(cf. equations \eqref{fpc:simple:3a}-\eqref{fpc:simple:3b}).
The KKT conditions provide a precise description of an optimal flow.
From now on the KKT conditions of \eqref{LP-TSO} will be called the
\emph{flow price condition}.

\begin{definition}[(Flow price condition)]\label{def:flow:price:condition}
Let $\tau\in\R^{C\times T}$, $\taulb\leq\tau\leq\tauub$, and $\p\in P^{A\times T}$. The
flow-price tuple $(\tau,\p)$ satisfies the \emphdef{flow price condition} if there exist 
$\muub,\mulb,\rhoub,\rholb\in\R^{C\times T}$ so that for all $c=\rs\in C$ and $t\in T$
\begin{align}
\muub\ct,\mulb\ct,\rhoub\ct,\rholb\ct&\geq 0\;,\label{efp:1}\\
(\muub\ct -\mulb\ct) +(\rhoub\ct -\rholb\ct) +[(\rholb\ctp -\rhoub\ctp)]+\prt-\pst&=0
\;,\label{efp:2}\\
(\tauct-\tauub)\;\muub\ct&=0\;,\label{efp:3}\\
(-\tauct+\taulb)\;\mulb\ct&=0\;,\label{efp:5}\\
(\tauct-\tau_{c,t-1}-\tauramp)\;\rhoub\ct&=0\;,\label{efp:4}\\
(\tau_{c,t-1}-\tauct-\tauramp)\;\rholb\ct&=0\;.\label{efp:6}
\end{align}
In terms of duality theory the conditions
\eqref{efp:1} and \eqref{efp:2} correspond to dual
  feasibility conditions and the remaining ones are
complementarity conditions. The term in square brackets is only added if $t+1\in T$.
\end{definition}

In the following we provide three implications that help to understand 
the properties of an optimal flow. A derivation of these implications
can be found in Appendix \ref{apx:fpc}.

If there is no ramping on an interconnector $c=\rs$, this condition
simplifies to
\begin{equation}
\forall\hours\qquad\pst-\prt\begin{cases}
	\leq 0&, \tauct<\tauub\\
	\geq 0&, \taulb<\tauct
\end{cases}\label{fpc:simple:1}.
\end{equation}
The economic interpretation of this simplified condition is that free transfer capacity in the
forward direction implies that the price $\prt$ in the source area must
be higher or equal to the price $\pst$ in the sink area. If this would not
be the case, then the free capacity should be used to transfer additional
electricity from the source to the sink, and therefore the market
would not be in an equilibrium.

If, on the other hand, the ramping condition (\ref{nb:ramprate}) on an interconnector
$c=\rs$ is only binding (i.e., satisfied with equality) in two
consecutive hours $t$ and $t+1$, and there is still available transfer
capacity in that time window, then the condition simplifies to:
\begin{equation}
\pst-\prt\enskip=\enskip-(\p_{s,t+1}-\p_{r,t+1})\;.\label{fpc:simple:2}
\end{equation}
This establishes that the price difference at one point in time can
have an impact on the price difference in the following ones.

As mentioned before the flow price condition implies that prices of adjacent areas can only
deviate if at least one respective transmission constraint is active
(terms in square brackets are only present if $t+1\in T$):
\begin{align}
&\forall\,c=\rs\in C,\:t\in T\notag\\
&\quad \prt<\pst\quad \Rightarrow\quad	\tau_{c,t}=\tauub\quad\vee\quad
	\tau_{c,t}-\tau_{c,t-1}=\tauramp\quad[\,\vee\quad
	\tau_{c,t}-\tau_{c,t+1}=\tauramp\,]\;,\label{fpc:simple:3a}
\\
&\quad \prt>\pst\quad \Rightarrow\quad	\tau_{c,t}=\taulb\quad\vee\quad
	\tau_{c,t-1}-\tau_{c,t}=\tauramp\quad[\,\vee\quad
	\tau_{c,t+1}-\tau_{c,t}=\tauramp\,]\;.\label{fpc:simple:3b}
\end{align}

\item{\textbf{Filling condition of linear bid curves.}}\label{nb:fillingcondition}
This constraint represents the optimality conditions of linear bid curves. Participants
who submit linear bid curves are modeled as price takers. In the optimization problem
associated with a linear bid curve $h$ the price is exogenously given and the only
remaining decision variable is $\delta_h$. The objective is to maximize the economic
surplus.
\begin{align*}
\max\quad&\int_0^{\delta_h}(\phvonq(u)-\pat)\,\dd u
\tag{QP-Hourly}\label{QP-Hourly}\\
\st \quad&[\vub]\quad+\delta_h\leq\dhub\\
         &[\vlb]\quad-\delta_h\leq\dhlb
\end{align*}
with $\phvond=\phd$. The variables in squared brackets correspond to the
dual variables of this problem.
From now on the KKT conditions of this price-parameterized optimization
problem will be called \emph{filling condition}.
\begin{definition}[(Filling condition)]\label{satz:Filling:Condition:kkt}
Let \deltawithinbounds\ and $\p\in P^{A\times T}$. The tuple $(\delta,\p)$
satisfies the \emphdef{filling condition} if and only if there exist
variables $\vub,\vlb\in\R^H$ so that for all \areas, \hours, $\hourly$
\begin{align}
	-\phvond+\pat+\vub_h-\vlb_h&=0\label{eq filling 1}\\
	(+\delta_h-\dhub)\;\vub_h &=0\\
	(-\delta_h+\dhlb)\;\vlb_h &=0\\
	\vub_h,\;\vlb_h&\geq 0\;.
\end{align}
\end{definition}
This condition can be reformulated to an equivalent condition:
\begin{proposition}[(Filling condition)]\label{satz:Filling:Condition}
Let \deltawithinbounds\ and $\p\in P^{A\times T}$. The tuple $(\delta,\p)$ satisfies the \emphdef{filling condition} if and only if
\begin{align}
&\forall\areas,\hours,\hourly\mit\dhlb<\dhub:
&\pat&\begin{cases}
 \leq\phvond, &\text{if } \delta_h=\dhub\\
 =\phvond, &\text{if }\delta_h\in (\dhlb,\dhub)\\
 \geq\phvond, &\text{if }\delta_h=\dhlb
\end{cases}\label{eq:FillingCond1}
\end{align}
with $\phvond=\phd$.
\end{proposition}
\begin{proof}
Let $a\in A,t\in T$, and $h\in\Hat$. If $\dhlb=\dhub$ then $\vlb_h$ and $\vub_h$
are free and \eqref{eq filling 1} is trivially satisfied.
Let $\dhlb<\dhub$. We will now consider three cases.
In the first case, let $\delta_h=\dhub$. Then $\vlb_h=0$ and
$\vub_h\geq0$, and \eqref{eq filling 1} simplifies to
$-\phvond+\pat\leq 0$. For the second case suppose
$\delta_h\in(\dhlb,\dhub)$. 
Then $\vlb_h=0$ and $\vub_h=0$, and \eqref{eq filling 1} simplifies to
$-\phvond+\pat=0$. For the last case, let $\delta_h=\dhlb$.
Then $\vlb_h\geq0$ and $\vub_h=0$, and \eqref{eq filling 1} simplifies
to $-\phvond+\pat\geq 0$.
\end{proof}
The second definition directly shows that an optimal solution only contains
quantity-price-combinations $(\delta,\p)$ that were defined by the piecewise
linear bid curves.
A MIP formulation of the filling condition based on the 
\emph{delta-method} (cf.~\cite[Chapter 26-3-I(f)]{dantzig1963}) can be found in
%\cite{markowitz1957}
Appendix~\ref{apx:filling:condition} and
a graphical illustration is given in Figure~\ref{fig:hourly:bids}.
\begin{figure}[hpb]
	\centering
	\begin{minipage}[c]{\textwidth}
\includegraphics[width=1.00\textwidth]{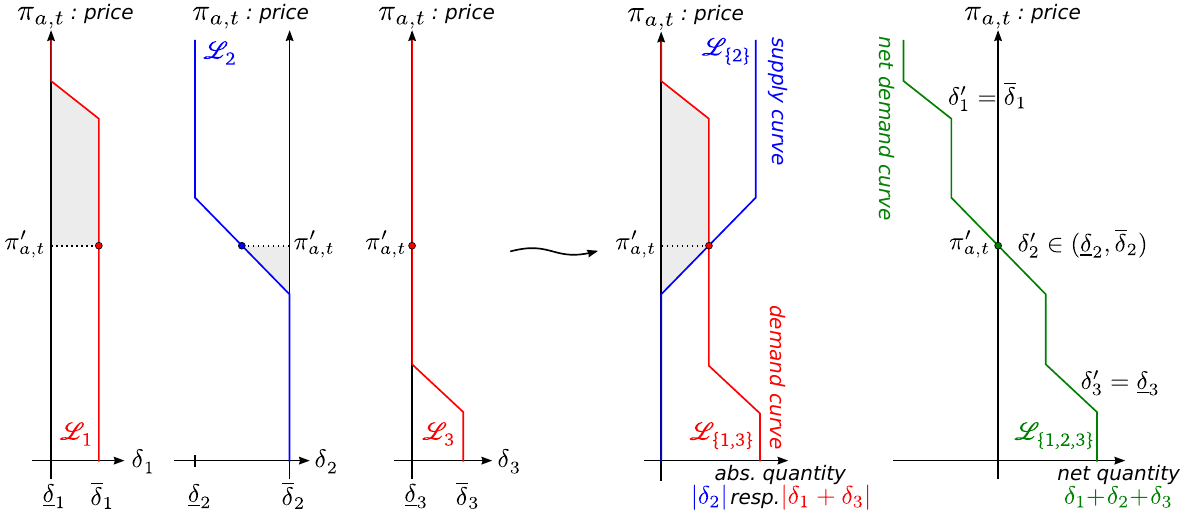}
	\centering
	\end{minipage}
\caption{%\vtop{\hbox{%
The tuple $(\delta',\p')$ satisfies the filling condition.
The points $(\delta'_h,\p'\at)$ are in $\mycal{L}_h$
for $h=1,2,3$ and the point $(\delta'_1+\delta'_2+\delta'_3,\p'\at)$
is in $\mycal{L}_{\{1,2,3\}}$.}
	\label{fig:hourly:bids}
\end{figure}

\end{enumerate}
\renewcommand{\theenumi}{\text{\arabic{enumi}}}
\renewcommand{\labelenumi}{\arabic{enumi}.}
\subsection{Additional market rules}
For the sake of completeness we need to mention that there exist additional
market rules that are not covered by the previous model. Some of these additional rules are only
valid for either the CWE region or the NPS region. For example in the CWE
region there is a parameter that prevents block bids with low marginal costs
or high marginal willingness to pay from being rejected.
In the NPS region participants can submit so called \emph{convertible block bids}.
A convertible block bid is basically a block bid,
but if it gets rejected, it is transformed into several
linear hourly bid curves and the auction starts again. Some interconnectors
are facing a significant transmission loss. The factor that models
this loss is called \emph{deadband}. There also exist rules
related to situations with extremely high or low area prices: 
The exchanges impose additional \emph{upper and lower bounds on the prices}.
If a price reaches a bound then \emph{curtailment rules} come into effect:
If for example an area price is at its upper bound, then demand block bids 
in this area have a lower priority than linear bid curves. 
Except for the different price intervals we will not explain these
rules in detail as they are out of the scope of this paper.

\begin{remark}\label{rem:priceinterval}
{\textit{Different price intervals.}}
The prices in the NPS region are limited to $[-200,2000]$ \eur/MW
whereas the prices in the CWE region are limited to $[-3000,3000]$ \eur/MW.
On the on hand the area prices \pat\ must be within the specific regional bounds.
On the other hand all submitted price limits (marginal willingness to pay
/ marginal cost) must be within the specific regional price interval:
Let $P_a$ be the feasible price interval in area $a$, then $p_b,p_f\in P_a$
for all block bids $b$ and flexible bids $f$ in area $a$, and 
$m_h([\dhlb,\dhub])\subseteq P_a$ for all linear bid curves $h$ in area $a$.
Note that linear bid curves are defined for all prices \pat\ in $\R$ and that they
are parallel to the price axis for prices in $\R\setminus P_a$.
%Figure~\ref{fig:priceintervals} shows an example.

If a solution contains a price $\pat\notin P_a$, then $\pat$ is changed
to the nearest price in $P_a$. By doing so, we might violate the flow
price condition. In other words, there might be a price
difference between two adjacent areas although none of the transport
restrictions is binding. In fact, this might happen if one of the areas has a small
price interval and reaches the maximal possible hourly demand or
supply.
\end{remark}

\subsection{Objective function}\label{section:objective}

The \surplus\ is the sum of consumer surplus and producer surplus
(cf. \cite{pindyck2005}).
It is also called \emph{global welfare}
(cf. \cite{kirschen2004}) or \emph{social welfare}.
For example a consumer buys $q$ units at a
price of $\p$ per unit and he is willing to pay a price of $p\ge\p$,
then his consumer surplus is given by $(p-\p)q$. The producer surplus
generated by hourly bid curves is determined accordingly by integrating 
the price per unit minus the producers marginal cost over the quantity
from zero to the produced quantity.

In Section \ref{sec:preliminaries} we already introduced the price-parameterized
surplus maximization problems of each trading product. We will now summarize
the results from above.% Let $\p\in\pricerange^{A\times T}$ be a linear price.

\begin{definition}
The \emphdef{surplus of a linear bid curve} $\hourly$ in area $a$ and hour $t$ is given by
\begin{align}
 \omega_h(\p,\delta)\label{eq:segment:surplus}
        :=&\int_0^{\delta_h}\phvonq(u)\,\dd u - \pat\delta_h
\end{align}
where $\phvond=\phd$.
Expression \eqref{eq:segment:surplus} is equivalent to the 
\emph{net consumers/producers surplus} in \cite[Chapter 2.2]{kirschen2004}.
Examples are provided in Figures~\ref{fig:willingness_to_pay}
and \ref{fig:marginal_cost}.
\end{definition}

\begin{definition}
The \emphdef{surplus of a block bid} \blocka\ with \areas\ for $\p\in
\pricerange^{A\times T}$ and $\beta\in\{0,1\}^B$ is given by
\[
 \omega_b(\p,\,\beta):=\sum_\hours p_b q\bt\beta_b - \sum_\hours \pat q\bt\beta_b.
\]
\end{definition}

\begin{definition}
The \emphdef{surplus of a flexible bid} $f\in F_a$ with $a\in A$ for
$\p\in\R^{A\times T}$ and $\varphi\in\{0,1\}^{F\times T}$ is given by
\[
 \omega_f(\p,\,\varphi):=\sum_\hours p_f q_f\pft - \sum_\hours\pat q_f\pft.
\]
\end{definition}

Another surplus arises on the interconnectors. This surplus represents the
\emphdef{surplus of cross-border trades} and is usually referred to as \emphdef{congestion rent}.

\begin{definition}\label{def:Rente:CBF}
The \emphdef{congestion rent} on interconnector $c=\rs\in C$ for $\p\in\R^{A\times T}\und\tau\in\R^{C\times T}$ is given by:
\[
 \omega\ct(\p,\,\tau):=\sum_\hours(\pst-\prt)\tau\ct\;.
\]
\end{definition}

Note that the congestion rent is indeed contained in the economic
surplus in most market models run today. The surplus generated due to congestion 
rent is usually reinvested in order to reduce tariffs or build new
lines\footnote{In 2007: Energinet.dk: 100\% to reduce tariffs, Fingrid: 100\%
  to build new lines, Statnett: 100\% to reduction, Svenska Kraftn\"at:
  100\% to build new lines.
(cf.~\href{http://www.nordpoolspot.com/Global/Download\%20Center/TSO/How-to-calculate-the-TSO-Congestion-rent.pdf}{http://www.nordpoolspot.com})}.
Later we will see that the flows are not uniquely determined by the
maximal \surplus{}, but in general ambiguities only exist, when prices of adjacent 
areas coincide. In these situations, there is no price difference ($\pst-\prt$) and 
thus no congestion rent, and the congestion rent cannot be influenced by changing 
the flow $\tau\ct$ within the ambiguities. For more details on ambiguities of flows 
and how to choose unique flows refer to Section~\ref{sec:choose-uniq-flows}.

\begin{definition}
Let \valP\ be feasible for \qpflow. The \emphdef{\surplus{} of the
  market} is the sum of the surpluses of all market participants and
is computed as follows:
\begin{align}
\label{eq:1}
\omega\varqpflow
=&\sum_\areas\sum_\hours\sum_\hourly\omega_h(\p,\delta)
    +\sum_\areas\sum_\blocka\omega_b(\p,\,\beta)
    +\sum_\areas\sum_\flexa\omega_f(\p,\,\varphi)
  \\\nonumber&
    +\sum_{c=\rs\in C}\omega\ct(\p,\,\tau)\\
=&
\label{eq:cost}
\sum_\areas\sum_\hours\left(
    \sum_\hourly\int_0^{\delta_h}\phvonq(z)\dd z
    + \sum_{\substack{\blocka}} p_b\,q_{b,t}\,\beta_b
    + \sum_{\substack{\flexa}} p_f\,q_f\,\pft\right)
  \\\nonumber&
    - \sum_\areas\sum_\hours\left(
    \sum_\hourly\pat\delta_h
    + \sum_{\substack{\blocka}} \pat q\bt\beta_b
    + \sum_{\substack{\flexa}} \pat q_f\pft\right)
  \\\nonumber&
    + \sum_{c=\rs\in C}\sum_\hours(\pst-\prt)\tau\ct
    \displaybreak[0]
\end{align}
\end{definition}
The first line of equation \eqref{eq:cost} represents the willingness to pay or
the cost of each bid, whereas the second and third represent the net amount of
money to be paid or received. The next theorem shows us that the latter is equal
to zero if the clearing condition \nbref{nb:clearing} is satisfied.
The simplified expression of the \surplus\ is equivalent to the objective of
the welfare maximization problem in \cite[Section II.B]{niu2005}.

\begin{theorem}\label{satz:Zielfunktion}
Let \valP\ be within the given bounds and satisfy condition
\nbref{nb:clearing}, then the \surplus\ is given by
\begin{align*}
\omega(\beta,\varphi,\delta)
&=\sum_\areas\sum_\hours\left(
 \sum_\hourly\int_0^{\delta_h}\phvonq(z)\dd z
 + \sum_{\substack{\blocka}} p_b\,q_{b,t}\,\beta_b
 + \sum_{\substack{\flexa}} p_f\,q_f\,\pft\right)\\
&=\sum_{h\in H}\Bigl(p_h\delta_h - \frac{1}{2}\dph\delta_h^2\Bigr)
 + \sum_{\substack{\blocks\\\hours}} p_b\,q_{b,t}\,\beta_b
 + \sum_{\substack{\flex\\\hours}} p_f\,q_f\,\pft\;.
\end{align*}
\end{theorem}
\begin{proof}
By construction we have
$C=\mathbin{\dot{\bigcup}}_\areas C^+_a=\mathbin{\dot{\bigcup}}_\areas C^-_a$.
With this equation we get for all \hours
\begin{equation}\label{lem:Kuppelleitungen:eq1}
\sum_{c=\rs\in C}(\pst - \prt)\tau_{c,t}
  = \sum_\areas \pat\Biggl(\sum_\cin\tau_{c,t} -
  \sum_\cout\tau_{c,t}\Biggr). 
\end{equation}
Combining the above with constraint \nbref{nb:clearing} and
substituting into \eqref{eq:1} %modified by Equation \eqref{eq:segment:surplus} 
yields the desired result.\hfill
\end{proof}

Observe that the objective function obtained in
Theorem~\ref{satz:Zielfunktion} is a quadratic, concave function and does not 
involve price variables $\p$  or flow variables $\tau$. While this is
convenient for the actual computation, these variables still have to
satisfy the stated conditions and so they implicitly influence the
\surplus{}.

The problem formulation \qpflow{}  with the objective function as
stated above is in principle solvable by convex Mixed-Integer Quadratic Programming
solvers like \cplex\ if we linearize the price conditions by using the big-M method. Solving
real world instances with this technique however is very hard, because
each non-convex quadratic price constraint must be modeled by using an auxiliary binary
variable.
For the considered instances \cplex\ was not able to obtain feasible solutions or improve 
provided warm-start solutions within 30 minutes.
We also observed that MINLP solvers and general MIQP solvers like 
\emph{BONMIN, COUENNE, SCIP}, and \emph{BARON} were not able to solve a
typical instance of the relaxed \qpbidcut\ model to optimality within 10 minutes%
\footnote{\emph{Run on a compute server with two 6-Core AMD Opteron 2435 (2.6GHz) CPUs, 64GB RAM, and 64-bit Debian}.}.
This indicates that these solvers will not be able to solve the the problem 
to optimality within 10 minutes if we add the non-convex complementarity conditions.
This is clearly unsatisfactory as the actual market coupling auction needs
to be cleared within 10 minutes. We therefore propose a heuristic and an 
exact algorithm for this problem.

%%%%%%%%%%%%%%%%%%%%%%
%%%% Optimality conditions
%%%%%%%%%%%%%%%%%%%%%%

\section{Optimality conditions}
\label{sec:optim-cond}
In this section we define a relaxation of \qpflow{} and analyze the
optimality conditions of this problem. Based on this analysis we
derive a heuristic and an exact algorithm that performs very well in practice.
We say that $\gebotswahl$ is a \emphdef{bid selection} if
$\valbeta\in\bin^B$ and $\valphi\in\bin^{F\times T}$. Note that a
bid selection does not necessarily satisfy the block price condition.

We consider the following parameterized and relaxed version of \qpflow{}.
We relax all price constraints and assume that a fixed bid selection \gebotswahl\
is exogenously given:
\begin{align*}
\max\qquad
		&%\omega(\delta,\tau)=
		\sum_{h\in H}\int_0^{\delta_h}\phvonq(u)\,\mathrm{d}u
			+\sum_{\blocks,\,\hours}p_b\,q_{b,t}\,\valbeta_b
			+\sum_{\flex,\,\hours} p_f\,q_f\,\valphi_{f,t}\\
\st\qquad
		&[\p]\quad\forall a\in A,t\in T
\qquad \sum_\hourly\delta_h+\sum_\cout\tauct
		-\sum_\cin\tauct
		+\sum_\blocka q_{b,t}\,\valbeta_b + \sum_\flexa q_f\,\valphi_{f,t}
		=0\\
		&[\vub]\quad\forall h\in H\hphantom{,c\in C}
\qquad\delta_h\leq\dhub\\
		&[\vlb]\quad\forall h\in H\hphantom{,c\in C}
\qquad-\delta_h\leq\dhlb\tag{QPRelax$_{\valbeta,\valphi}$}\label{QPRelax}\\
		&[\muub]\quad\forall c\in C,t\in T
\qquad \tauct-\tauub\leq0\\
		&[\rhoub]\quad\forall c\in C,t\in T
\qquad \tauct-\tau_{c,t-1}-\tauramp\leq0\\
		&[\mulb]\quad\forall c\in C,t\in T
\qquad -\tauct+\taulb\leq0\\
		&[\rholb]\quad\forall c\in C,t\in T
\qquad \tau_{c,t-1}-\tauct-\tauramp\leq0
\end{align*}

Note that \qprelax{} is a parameterized optimization problem where
the binary variables $\gebotswahl$ are exogenously given parameters.
It is a convex optimization problem with a polyhedral feasible region.
Therefore, it is well defined to associate dual variables with every
primal constraint. The dual variables are denoted by the Greek letters
that are given in squared brackets.
The next theorem will show that in this relaxation we do
not need to require the filling condition \nbref{nb:fillingcondition}
or the flow price condition \nbref{nb:fpc} explicitly as these
naturally follow from the KKT conditions.

\begin{theorem}[(Natural spot price characteristics)]\label{satz:spotpreiseigenschaften}
Let $\gebotswahl$ be an exogenously given fixed bid selection, i.e., $\valbeta$ and $\valphi$ are
exogenously given parameters. A feasible solution $(\delta,\tau)$ to \qprelax\
is an optimal solution if and only if there exist prices $\p$ so that $(\delta,\p)$
satisfies the filling condition and $(\tau,\p)$ satisfies the flow price condition.
These prices are given by the optimal dual variables to the clearing condition of
\qprelax.
\end{theorem}

\begin{proof}
Let \gebotswahl\ be a fixed bid selection. The objective of the parameterized
optimization problem \qprelax\ is concave and continuous differentiable and 
the constraints are affine linear. The Karush-Kuhn-Tucker optimality 
conditions provide the following: A primal feasible solution $(\delta,\tau)$
to \qprelax\ is an optimal solution if and only if there exist dual 
variables $\p,\vub,\vlb,\muub,\mulb,\rhoub,\rholb$ with
\begin{align}
&&\vub,\vlb,\muub,\mulb,\rhoub,\rholb&\geq0
		\label{eq:natural:1}\\
&\forall\areas,\hours,\hourly
	&-\phvonq(\delta_h)+\pat+\vub_h-\vlb_h&=0
		\label{eq:natural:2}\\
&\forall c=\rs\in C,\hours
&(\muub\ct -\mulb\ct) +(\rhoub\ct -\rholb\ct) +[(\rholb\ctp -\rhoub\ctp)]&=\pst-\prt
			\label{eq:natural:3}\\
&\forall\areas,\hours,\hourly
	&(\delta_h-\dhub)\;\vub_h &= 0\label{eq:natural:4}\\
	&&(-\delta_h+\dhlb)\;\vlb_h &= 0\label{eq:natural:4:ub}\\
&\forall c=\rs\in C,\hours
	&(\tauct-\tauub)\;\muub\ct &= 0\label{eq:natural:5}\\
	&&(\tauct-\tau_{c,t-1}-\tauramp)\;\rhoub\ct &= 0\\
	&&(-\tauct+\taulb)\;\mulb\ct &= 0\\
	&&(\tau_{c,t-1}-\tauct-\tauramp)\;\rholb\ct &= 0,\label{eq:natural:5:ub}
\end{align}
where the term in squared brackets is only added if \(t+1 \in T\).
The conditions \eqref{eq:natural:1} to \eqref{eq:natural:3} 
are the dual feasibility conditions and conditions
\eqref{eq:natural:4} to \eqref{eq:natural:5:ub} are the
complementarity conditions. Observe that equations \eqref{eq:natural:1}, \eqref{eq:natural:3}, and
\eqref{eq:natural:5} to \eqref{eq:natural:5:ub} imply that $(\tau,\p)$
fulfills the flow price condition
\nbref{nb:fpc}. The filling condition
\nbref{nb:fillingcondition} corresponds to the equations
\eqref{eq:natural:1}, \eqref{eq:natural:2}, and \eqref{eq:natural:4}
to \eqref{eq:natural:4:ub}.
\end{proof}

We will later use Theorem~\ref{satz:spotpreiseigenschaften} to extend solutions to
\qprelax{} to feasible ones for \qpflow{}.

%\input{4_optimality_conditions}

%%%%%%%%%%%%%%%%%%%%%%
%%%% Uniqueness
%%%%%%%%%%%%%%%%%%%%%%

\section{Uniqueness of a solution}
\label{sec:uniqueness-solution}
In the following section we will discuss uniqueness aspects of the
computed solution. More precisely, we will show that the economic
surplus for a given bid selection is unique. The flows and the prices however are not
necessarily unique as we will see. Suppose that an algorithm finds
an optimal solution to (MPEC). We will see that if we fix the binary variables
in (MPEC) to the optimal values, then all feasible solutions to the remaining model
have the same objective value, thus the remaining model is just a feasibility problem.
In practice there exist several feasible solutions to the remaining model.
For this reason an auctioneer must define rules that determine a unique solution
among these feasible solutions. In practice SESAM and COSMOS use
slightly different rules to choose unique solutions. For the sake of
exposition we propose a simple rule that chooses unique flows and prices without
changing the total economic surplus. At first an algorithm chooses a bid selection,
then a unique flow is chosen, and finally a unique price.
As a corollary of \ref{satz:spotpreiseigenschaften} we obtain:

\begin{corollary}\label{cor:qprelax:optimality}
Let \gebotswahl\ be a fixed bid selection and let \valrelax\ be a feasible
solution to the parameterized optimization problem \qprelax. Then \valrelax\
is an optimal solution to \qprelax\ if and only if there exist prices $\valp$ 
that satisfy the filling condition and the flow price condition.
\end{corollary}

\subsection{Unique economic surplus for a given bid selection}
Let \gebotswahl\ be a bid selection. 
Recall that the model \qpflow{} requires the filling condition and the
flow price condition to be satisfied. We can therefore only construct a solution to
\qpflow{} from a solution to \qprelax{} if there exist prices that satisfy the filling and the flow
price condition. Fortunately, these are exactly the KKT optimality conditions of
\qprelax. Let \valrelax\ be a feasible solution to \qprelax\ and suppose
that there exist prices that satisfy the filling condition and the
flow price condition. Corollary \ref{cor:qprelax:optimality} provides
that in this case \valrelax\ is an optimal solution to \qprelax.
There exists exactly one optimal objective value for each parameterized
\qprelax. If we fix the bid selection in the model \qpflow\
then the remaining degrees of freedom cannot influence the \surplus\ anymore.

\subsection{Uniqueness of hourly bid execution}
Let \gebotswahl\ be a bid selection and let $I\subset H$ be the indices
of linear bid curves with a strict decreasing marginal willingness to
pay / costs function $m_h$, i.e., 
\[
I:=\set{h\in H\mid \dph\neq 0}\;.
\]
%Note that $H$ only contains non-horizontal segments, i.e., $\dhlb<\dhub$.
An optimal solution to \qprelax{} is unique with respect to $\delta_I$,
as $\omega$ is strictly concave with respect to $\delta_I$ (cf.~\cite{dattorro05}).

\subsection{Choosing unique flows}
\label{sec:choose-uniq-flows}
\begin{figure}[tb]
	\centering
	\begin{minipage}[c]{10cm}
	\includegraphics[width=1.00\textwidth]{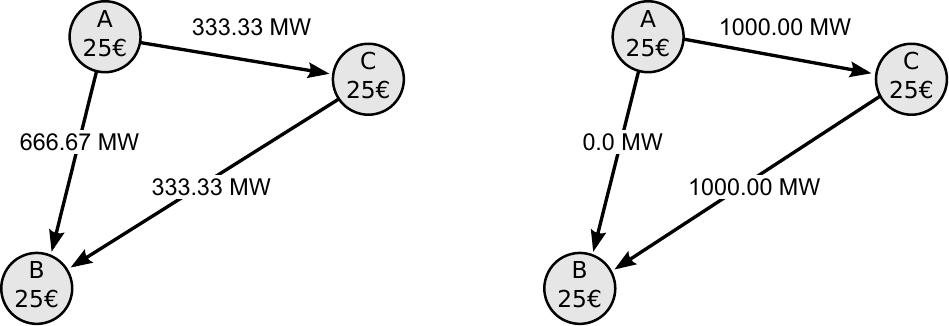}
	\centering
	\end{minipage}
	\caption{Both flows generate the same \surplus.}
	\label{fig:Eindeutigkeit}
\end{figure}

The optimality of a solution to \qprelax\ and the 
constraints in Section \ref{section:constraints} are not sufficient
to imply a unique flow. This follows from the non-uniqueness of the execution
state of some linear bid curves (those with \(\dph = 0\)).
To obtain a unique flow, we minimize the squared
flow while fixing the \surplus. In Figure~\ref{fig:Eindeutigkeit} an
example of a non-unique flow is depicted. The left part shows the
solution with minimized squared flows. In this example the prices
of the adjacent areas coincide and there is no congestion rent.

The following model selects a unique flow: Let \gebotswahl\ be a
bid selection, \valrelax\ be an optimal solution to
\qprelax{}, and let $J\subset\set{h\in H\mid \dph=0}$
be the set of linear bid curves with non-unique execution state.

\begin{align*}
\min\quad&\sum_\hours\sum_\cons\tauct^2\tag{FixFlow}\label{FixFlow}\\
\st\quad&\text{The \surplus\ must not change (cf.~Theorem \ref{satz:Zielfunktion}):}\\
&\quad\sum_{h\in J}p_h\delta_h
=\sum_{h\in J}p_h\valdelta_h\\
&\text{Clearing condition (cf. \nbref{nb:clearing}):}\\
&\quad\forall\areas,\hours\\
&\quad\qquad\sum_{h\in J\cap\Hat}\delta_h 
+\sum_\cout\tauct-\sum_\cin\tauct
=\sum_{h\in J\cap\Hat}\valdelta_h 
+\sum_\cout\valtauct-\sum_\cin\valtauct\\
&\text{Upper bounds, lower bounds, and ramp rate:}\\
&\quad\forall h\in J\qquad\quad\qquad \dhlb\leq\delta_h\leq\dhub\\
&\quad\forall\cons,\hours\qquad \taulb\leq\tauct\leq\tauub\\
&\quad\forall\cons,\hours\qquad -\tauramp\leq\tauct-\tau\ctm\leq\tauramp
\end{align*}

Note that the objective is strictly convex and the constraints are
affine linear. So the obtained solution is unique with respect to the flow.
The unique flow determines a unique net export per area and hour, therefore
also the execution state of the linear bid curves in $J$ is unique.

\subsection{Choosing unique prices}
In this section we will explain a way of how to obtain unique prices. Let
\gebotswahl\ be a bid selection and \valrelax\ an optimal solution
to \qprelax\ with unique flows. Then there exist prices $\valp$ that
satisfy the filling condition and the flow price condition. The following
example will illustrate that these prices are not unique in general.

Suppose that there is only one market area $a$, and one time slot $t$.
Let the hourly bid curve $\mycal{L}_{H\at}$ be parallel to the 
price axis in a price interval $I$ (i.e., no volume change).
Let $\valdelta$ be an optimal solution to \qprelax\ and let
$\valp$ be a price that satisfies the filling condition. The
flow price condition is trivially satisfied as there exists no
interconnector. If $\valp$ is in $I$ then all prices in $I$
satisfy the filling condition.
%Figure \ref{fig:nonunique:price} illustrates this example.

The following model takes as input a given bid selection and
an optimal solution to \qprelax, and then determines unique prices
by minimizing the squared prices subject to all price conditions.

\begin{align*}
 \min\enskip
	&\text{\it (squared prices)}
	&&\sum_\areas\sum_\hours \pat^2
		\tag{QPPrice}\label{QPPrice}\\
\st\enskip
&\text{\it (price constraints)}\\
&\forall\areas, \blocka&&
		\valbeta_b\sum_\hours(p_b - \pat)\,q_{b,t}\geq 0
		\tag{\ref{nb:blockpreisbedingung}}\\
&\forall\areas, \flexa, \hours&&
		\valphi_{f,t}(p_f-\pat)\,\sgn(q_f)\geq 0
		\tag{\ref{nb:flexpreisbedingung}}\\
&\forall c=\rs\in C,\hours&&
%		(\muub\ct -\mulb\ct) +(\rhoub\ct -\rholb\ct) +[(\rholb\ctp -\rhoub\ctp)]+\prt-\pst=0
		(\muub\ct -\mulb\ct) +(\rhoub\ct -\rholb\ct) +[(\rholb\ctp -\rhoub\ctp)]=\pst-\prt
		\tag{\ref{nb:fpc}}\\
&\forall\cons,\hours&&
		(\valtauct-\tauub)\;\muub\ct=0
		\tag{\ref{nb:fpc}}\\
&\forall\cons,\hours&&
		(-\valtauct+\taulb)\;\mulb\ct=0
		\tag{\ref{nb:fpc}}\\
&\forall\cons,\hours&&
		(\valtauct-\valtau_{c,t-1}-\tauramp)\;\rhoub\ct=0
		\tag{\ref{nb:fpc}}\\
&\forall\cons,\hours&&
		(\valtau_{c,t-1}-\valtauct-\tauramp)\;\rholb\ct=0
		\tag{\ref{nb:fpc}}\\
&\forall\cons,\hours&&
		\muub\ct,\;\mulb\ct,\;\rhoub\ct,\;\rholb\ct\geq 0
		\tag{\ref{nb:fpc}}\\
&\forall\areas,\hours,\hourly&&
		-\phvonq(\valdelta_h)+\pat+\vub_h-\vlb_h=0
		\tag{\ref{nb:fillingcondition}}\\
&\forall\areas,\hours,\hourly&&
		(\valdelta_h-\dhub)\;\vub_h = 0
		\tag{\ref{nb:fillingcondition}}\\
&\forall\areas,\hours,\hourly&&
		(-\valdelta_h+\dhlb)\;\vlb_h = 0
		\tag{\ref{nb:fillingcondition}}\\
&\forall\areas,\hours,\hourly&&
		\vub_h,\;\vlb_h\geq 0
		\tag{\ref{nb:fillingcondition}}\\
&\forall\areas,\hours&&\pat\in\pricerange
\end{align*}

Observe that the filling condition, the flow price condition, and the
bid price constraints simplify to linear constraints, because the flow
$\valtau$, the hourly bid execution $\valdelta$, and the bid selection
$\gebotswahl$ are fixed.

%%%%%%%%%%%%%%%%%%%%%%
%%%% Algorithm
%%%%%%%%%%%%%%%%%%%%%%

\section{The bid cut heuristic and an exact algorithm}
\label{sec:algorithm}
We will now present the proposed algorithm to solve the clearing
problem. The basic idea is to relax all price
constraints. For example we allow that block or flex
bids incur a loss. The following model is similar to \qprelax,
but this time the bid selection is not fixed: The model can choose
the best bid selection in terms of the total \surplus.
We will see that an optimal solution to this model has similar
properties to those in Theorem \ref{satz:spotpreiseigenschaften}.
\begin{align*}
 \max\enskip
	&\text{\it (\surplus)}
	&&\sum_{h\in H}\int_0^{\delta_h}\phvonq(u) \dd u
		+\sum_{\substack{\blocks\\\hours}} p_b\,q_{b,t}\,\beta_b
		+\sum_{\substack{\flex\\\hours}} p_f\,q_f\,\pft
		\label{QPBidCut}\tag{QPBidCut}\\
\st\enskip
&\text{\it (quantity constraints)}\\
&\forall a\in A,\,t\in T	&&
		\sum_\hourly\delta_h +\sum_\blocka q_{b,t}\beta_b
		 + \sum_\flexa q_f \pft
		=\sum_\cin \tau_{c,t} - \sum_\cout \tau_{c,t}
		\tag{\ref{nb:clearing}}\\
\enskip&\forall c\in C,\,t\in T	&&
 		-\tauramp\leq\tau_{c,t} - \tau_{c,t-1}\leq\tauramp
 		\tag{\ref{nb:ramprate}}\\
&\forall (b,c)\in L&&	\beta_b\leq\beta_c
		\tag{\ref{nb:blocklinks}}\\
&\forall f\in F&&	\sum_\hours \pft\leq 1
		\tag{\ref{nb:flexphi}}\\
&\forall\blocks&&	\beta_b\in\{0,1\}\\
&\forall\flex,\hours&&	\pft\in\{0,1\}\\
&\forall\cons,\hours&&	\tau_{c,t}\in[\taulb,\tauub]\\
&\forall h\in H&&	\delta_h\in[\dhlb,\dhub]\\
\end{align*}
\begin{theorem}\label{satz:BidCutEigenschaften}
Let $\valbidcut$ be an optimum solution to \eqref{QPBidCut}.
Then there exist prices $\p$ that satisfy the filling condition
and the flow price condition.
The objective value $\omega(\valbeta,\valphi,\valdelta)$ is
the \surplus\ in the point $\valbidcut$.
\begin{proof}
Let \valbidcut\ be optimal for \eqref{QPBidCut}. Then
$(\valbeta,\valphi)$ is a fixed bid selection and \valrelax\ is
optimal for \eqref{QPRelax}. Then Theorem \ref{satz:spotpreiseigenschaften}
yields the existence of prices that satisfy the filling condition and
the flow price condition. Every feasible solution to \eqref{QPBidCut} 
satisfies the clearing condition \eqref{nb:clearing}. Therefore, Theorem
\ref{satz:Zielfunktion} yields that $\omega(\valbeta,\valphi,\valdelta)$
is the economic surplus.
\end{proof}
\end{theorem}

Note that the theorem above only guarantees that the filling condition
and the flow price condition can be satisfied in an optimal solution.
It does not guarantee that none of the executed block bids or flexible bids
incurs a loss. 

Suppose now that we have an optimal solution to \qpbidcut. Then we can find
prices that satisfy the filling condition and the flow %\eqref{nb:fillingcondition}
price condition, but there might exist block or flexible bids that %\eqref{nb:fpc}
incur a loss. Thus the solution is not feasible for \qpflow. In this case
an additional constraint, a so called \emph{bid cut}, is added to
the problem. This constraint cuts off the infeasible
bid selection and a class of similar bid selections. This process is
repeated until a selection is found that is feasible for
\qpflow{}. The feasibility of a selection can be checked with a model that is
similar to \qpprice. As shown in Theorem \ref{satz:BidCutEigenschaften} 
the flow price condition \eqref{nb:fpc} and the filling condition \eqref{nb:fillingcondition}
can be satisfied. The block and flexible bid price condition must be relaxed.
Instead of choosing unique prices, we are now interested in minimizing
the violation of these two conditions. Model \eqref{LPPrice} shows these modifications.
We introduced a positive continuous slack variable $\lambda_b$ for every 
executed block or flexible bid $b$ and used them to relax the
according block or flexible bid price condition. Instead of minimizing
the squared prices we minimize the sum of these slack variables.
Now the variable $\lambda_b$ models the loss of the executed block or flexible bid $b$.

Let \valbidcut\ be an optimal solution to \qpbidcut, then the parameterized
linear program \eqref{LPPrice} checks whether the solution is feasible to \qpflow:
If there exists a solution with $\lambda=0$, then we have found a feasible solution
to \qpflow.
\begin{align*}
 \min\enskip
	&\text{\it (incurred loss)}
	&&\sum_\blocks \lambda_b +\sum_\flex \lambda_f
		\tag{LPPrice}\label{LPPrice}\\
\st\enskip
&\text{\it (price constraints)}\\
&\forall\areas, \blocka&&
		\valbeta_b\sum_\hours(p_b - \pat)\,q_{b,t}\geq -\lambda_b
		\tag{relaxed \ref{nb:blockpreisbedingung}}\\
&\forall\areas, \blocka&&\lambda_b\ge 0\\
&\forall\areas, \flexa, \hours&&
		\valphi_{f,t}(p_f-\pat)\,\sgn(q_f)\geq -\lambda_f
		\tag{relaxed \ref{nb:flexpreisbedingung}}\\
&\forall\areas, \flexa&&\lambda_f\ge 0\\
&\forall c=\rs\in C,\hours&&
%		(\muub\ct -\mulb\ct) +(\rhoub\ct -\rholb\ct) +[(\rholb\ctp -\rhoub\ctp)]+\prt-\pst=0
		(\muub\ct -\mulb\ct) +(\rhoub\ct -\rholb\ct) +[(\rholb\ctp -\rhoub\ctp)]=\pst-\prt
		\tag{\ref{nb:fpc}}\\
&\forall\cons,\hours&&
		(\valtauct-\tauub)\;\muub\ct=0
		\tag{\ref{nb:fpc}}\\
&\forall\cons,\hours&&
		(-\valtauct+\taulb)\;\mulb\ct=0
		\tag{\ref{nb:fpc}}\\
&\forall\cons,\hours&&
		(\valtauct-\valtau_{c,t-1}-\tauramp)\;\rhoub\ct=0
		\tag{\ref{nb:fpc}}\\
&\forall\cons,\hours&&
		(\valtau_{c,t-1}-\valtauct-\tauramp)\;\rholb\ct=0
		\tag{\ref{nb:fpc}}\\
&\forall\cons,\hours&&
		\muub\ct,\;\mulb\ct,\;\rhoub\ct,\;\rholb\ct\geq 0
		\tag{\ref{nb:fpc}}\\
&\forall\areas,\hours,\hourly&&
		-\phvonq(\valdelta_h)+\pat+\vub_h-\vlb_h=0
		\tag{\ref{nb:fillingcondition}}\\
&\forall\areas,\hours,\hourly&&
		(\valdelta_h-\dhub)\;\vub_h = 0
		\tag{\ref{nb:fillingcondition}}\\
&\forall\areas,\hours,\hourly&&
		(-\valdelta_h+\dhlb)\;\vlb_h = 0
		\tag{\ref{nb:fillingcondition}}\\
&\forall\areas,\hours,\hourly&&
		\vub_h,\;\vlb_h\geq 0
		\tag{\ref{nb:fillingcondition}}\\
&\forall\areas,\hours&&\pat\in\pricerange
\end{align*}

We will now introduce two families of cuts that can be used. The first
one works very well in practice but is inexact as it might
converge to a suboptimal solution. The second family of cuts
is exact but more cuts have to be added. We will later report timings
for both families in Section~\ref{sec:results}. 

\begin{definition}
Let \valbidcut\ be an optimal of \qpbidcut\ and let $\valp$ be an optimal solution
to \eqref{LPPrice}. The set of bids that incur a loss at prices $\valp$ is given by
\begin{align*}
\Bloss
%&=\{\blocks\mid b\text{ is executed and incurs a loss at prices $\valp$}\}\\
&=\{\blocks\mid\valbeta_b=1\und\valblockgain < 0\}\\
\text{and }\Floss
%&=\{\flex\mid f\text{ is executed and incurs a loss at prices $\valp$}\}\\
&=\{\flex\mid\text{ex. }\hours\mit\valphi_{f,t}=1\und\valflexgain < 0\}.
\end{align*}
where \(\Bloss\) is the set of executed block bids that incur a loss at
prices \(\valp\) and \(\Floss\) is the set of executed flexible bids that incur a loss at
prices \(\valp\). The linear constraint that prohibits the execution
of at least one bid of the set $\Bloss\cup\Floss$ is given by
\[
\cut(\Bloss,\Floss)
:\qquad\sum_{b\in\Bloss}\beta_b + \sum_{f\in\Floss,t\in T:\valphi_{f,t}=1}\pft\leq |\Bloss|+|\Floss|-1.
\]
We refer to this class of constraints as \emphdef{Bid Cuts} (cf. \cite[Definition 5.2]{muel2009}).
\end{definition}

As mentioned above, for this family of cuts, the obtained solution
might be slightly suboptimal, i.e., these cuts provide us with a
heuristic. Nonetheless, this cut can be applied iteratively as
shown in Algorithm \ref{alg:Einfache:BidCutHeuristik} or it can be incorporated
into a branch-and-cut framework. The latter is the more effective one
as it provides a very fast heuristic: The solving process of \cplex{}
is only started once and the bid cuts are injected directly into the 
process by using callback functions.

If the Bid Cut is replaced with a less aggressive one that only removes
exactly one invalid bid selection at a time, then the resulting algorithm will
converge to the globally optimal solution. 
For example, the following family of cuts separates
exactly one invalid bid selection. 

\begin{definition}\label{def:no:good}
Let \valbidcut\ be a solution of \qpbidcut{}. Then following inequality cuts off the bid
selection \((\valbeta,\valphi)\):
\[
\cut(\valbeta,\valphi)
:\qquad\sum_{b\in B: \valbeta_b =0}\beta_b + \sum_{f \in F, t\in
  T:\valphi_{f,t}=0} \pft + \sum_{b\in B: \valbeta_b =1} (1-\beta_b) + \sum_{f \in F, t\in
  T:\valphi_{f,t}=1} (1-\pft) \geq 1
%\varphi_f \sum_{f\in\Floss,t\in T:\valphi_{f,t}=1}\pft\leq |\Bloss|+|\Floss|-1.
\]
We refer to this class of constraints as \emphdef{exact bid cuts}
(cf. \cite[Theorem 5.4]{muel2009}, ``no-good cuts'' in \cite{chu2004}).
\end{definition}

The most effective way to implement an exact algorithm that uses the exact bid cuts is to
inject them into a branch-and-cut framework. We decided to call this algorithm
\emph{Branch-And-Cut Decomposition}, because we decompose the model into the upper level
problem \qpbidcut\ and the parameterized lower level problems \eqref{LPPrice} and use
exact cuts to connect both levels.

\begin{algorithm}
\begin{algorithmic}
\medskip
\REQUIRE An instance of the model \qpbidcut
\ENSURE \valqpflow\ feasible for \qpflow
\STATE
\STATE $done \gets \FALSE$
\STATE $i\gets 0$
    \STATE
\WHILE{$\neg done$}
    \STATE
   \STATE$\valbidcut~\gets$ solve \qpbidcut.
   \STATE$(\valtau)~\gets$ solve (FixFlow) with input values $(\valbeta,\valphi,\valdelta)$.
   \STATE$(\valp,\Blossi,\Flossi)~\gets$ solve \eqref{LPPrice}\ with input values $\valbidcut$.
    \STATE
    \IF {$|\Blossi\cup\Flossi| = 0$,}
        \STATE $done\gets\TRUE$
    \ELSE
	\STATE add $\text{Cut}(\Blossi,\Flossi)$ to the model \eqref{QPBidCut}.
    \ENDIF 
    \STATE
    \STATE $i \gets i+1$
\ENDWHILE
\medskip
\STATE {\it Return} $(\valp,\valbeta,\valphi,\valdelta,\valtau)$.
\medskip
\end{algorithmic}
\caption{Iterative Bid Cut Heuristic}
\label{alg:Einfache:BidCutHeuristik}
\end{algorithm}

%\input{6_algorithm}

%%%%%%%%%%%%%%%%%%%%%%
%%%% Results
%%%%%%%%%%%%%%%%%%%%%%

\section{Results}
\label{sec:results}
We will now present some computational results for the algorithms
presented in Section~\ref{sec:algorithm}. We compare the Bid Cut
Heuristic and Branch-and-Cut Decomposition with the commercial algorithm 
\emph{EMCC Optimizer} of Deutsche B\"orse Systems which is currently used by 
European Market Coupling Company (EMCC) to determine
flows between the NPS and the CWE region.

For our computational tests we used 79 realistic instances that contain 10 European
market areas, about 600 combinatorial bids, and about $31700$ linear bid curves.
All tests were performed on the same hardware\footnote{Intel
  Xeon Core E5440 with 8 GB memory} using \cplex\ 12.1 as mixed
integer linear/quadratic programming solver. In our tests, we did not incorporate the 
\fixflow{} model, as its sole purpose is the redistribution of
flow without affecting the overall \surplus{}. Moreover, our tests
indicated that fixing the flow impacts the prices only in very
pathological cases that can be disregarded in this paper. We will first give
a brief introduction to an early version of the EMCC Optimizer
and provide some statistics. Then we discuss the two proposed
algorithms from above.

\subsection{EMCC Optimizer}
The version of the EMCC Optimizer that we used in our comparison is version 2.1.2
that was used in 2009. This algorithm is based on computing a start
solution with a linear program (cf.~\cite{krion2008})
and then checking the feasibility of the bid selection. If the bid selection 
is not feasible, executed bids are excluded successively until the 
solution is feasible. Then the solution will be improved by trying to 
include not executed bids. For more details see 
\cite[Chapter~4.3.2]{emccopt}. In some cases this search phase is very
time consuming and cannot be completed within 10 minutes\footnote{The
limit of 10 minutes reflected the operational requirements when the tests of this
algorithm were performed (September 2009). For production it has been changed to 15 minutes.}.
In this case, the last feasible solution will be returned; in fact the start solution
for the last phase is already feasible however might be suboptimal. In our analysis, for 22\% of the
instances this last phase could not be completed on time. The average time 
for computing the solutions amounts to 6 minutes and in at least 15\% of the cases the bid 
selection was provably optimal. In 76\% of the cases the selections could be slightly 
improved by our algorithms and in the remaining cases we could neither
prove nor disprove the optimality of the selection.
The relative distance to 
the best upper bound, i.e., the relative gap, averages to 2.577\e{-6};
a selection is considered optimal here if the relative gap is below 1\e{-12}
% actually it corresponds to $0.1$ euro cent
which roughly corresponds to one euro cent.
On average there are 10.76
rejected combinatorial bids, so called \emph{paradoxically rejected bids (PRB)}
(cf. \cite{meeus2009}) that could potentially have generated a net
profit. However, in most cases these bids cannot be included without
decreasing the overall \surplus\ or rejecting other bids that are in
the money. 

Some of our ideas where used to improve later versions of the EMCC solver and, 
therefore, increased the daily economic surplus indirectly. For example
the usage of only one aggregated hourly bid curve per area and hour instead
of separate demand and supply curves decreased the computing time significantly
and allowed for testing more solution candidates in less time.

\subsection{Bid Cut Heuristic}

The tests are based on a branch-and-cut version of the Bid Cut 
Heuristic. In at least 38\% of the tests the absolute gap is smaller than one
euro cent. The average relative gap amounts to 1.926\e{-6}. In 96\% of the cases, the solutions
found by the Bid Cut Heuristic could not be outperformed by any
other algorithm of the test. On average there are only 4.65 PRBs leading
to a potentially higher overall acceptance of the prices and executed quantities.
The most important property of the algorithm however is its running
time which averages to only 4.1 seconds. The maximum 
computing time was 1.1 minute and so clearing within the 10 minutes is
easily possible. Moreover, due to being able to compute the solution
way before the time limit is reached the auction can also be run on
different machines with identical results. This would not be possible
if the process would be stopped by the time limit as in this case the
actual performance of the machine would potentially impact the
computed prices. Note that reproduction of prices is crucial for market
transparency and confidence and an algorithm that is terminated
by a time limit is less suited for market clearing as the results
can be influenced by side effects like CPU load that is generated by other processes.

\subsection{B\&C Decomposition}

The Branch-and-Cut Decomposition is slower than the Bid Cut
Heuristic. However using it, we can prove that the absolute gap of a solution is smaller
than one euro cent. It uses the solutions of the Bid Cut Heuristic as warm start
solutions. In 38\% of the cases the algorithm finishes in about 9
seconds by finding the optimal solution and proving optimality. In the 
other cases the solver was stopped when reaching the time limit without reducing
the absolute gap to one euro cent. During the computation 
4\% of the bid selections found by the Bid Cut Heuristic could be
improved. This reduced the relative gap to 1.924\e{-6} and the number of PRBs 
was reduced to 4.62. The improvement which is not significant
indicates that the solutions obtained via Bid Cut 
Heuristic are very good.

%\subsection*{Relax And Fix}
%For the sake of completeness we also listed the results of our fastest heuristic,
%the \emph{Relax And Fix} heuristic. However the details of this heuristic are out
%of the scope of the paper, because the Bid Cut Heuristic determines solutions with
%a better \surplus.\\\medskip

\begin{table}[tpb]
\makebox[\textwidth][c]{
\begin{threeparttable}
	\centering{\Times{7.5}
	\begin{tabular}{|l|r|r|r|r|r|r|r|r|}
\hline
\multicolumn{ 1}{|c|}{\textit{Algorithm}} & \multicolumn{ 1}{c|}{\textit{Abs. Gap}} & \multicolumn{ 1}{c|}{\textit{Best known}} & \multicolumn{ 1}{c|}{\textit{PRBs}} & \multicolumn{ 1}{c|}{\textit{Relative}} & \multicolumn{ 3}{c|}{\textit{Computing Time $\le$ 10 Min}} & \multicolumn{ 1}{c|}{\textit{Stopped at}} \\ \cline{ 6- 8}
\multicolumn{ 1}{|c|}{} & \multicolumn{ 1}{c|}{\textit{$\le$ 0.01 EUR}} & \multicolumn{ 1}{c|}{\textit{Bid Selection}} & \multicolumn{ 1}{c|}{\textit{}} & \multicolumn{ 1}{c|}{\textit{Gap}} & \multicolumn{1}{c|}{\textit{Average (s)}} & \multicolumn{1}{c|}{\textit{Min (ms)}} & \multicolumn{1}{c|}{\textit{Max (s)}} & \multicolumn{ 1}{c|}{\textit{max Time}} \\ \hline\hline
EMCC Optim. 2.1.2 & $\geq$ 15\% & 24\% & 10.76 & 2.577E-6 & 361.4 & 46000 & 600.0 & 22\% \\ \hline
%Relax And Fix\tnote{1} & 
%$\geq$ 35\% & 71\% & 5.04 & 2.013E-6 & 0.7* & 193* & 1.6* & 0\% \\ \hline
Bid Cut Heuristic & $\geq$ 38\% & 96\% & 4.65 & 1.926E-6 & 4.1 & 493 & 62.8 & 0\% \\ \hline
B\&C Decomposition & $\geq$ 38\% & 100\% & 4.62 & 1.924E-6 &  (9.3) 390.8 & 1114 & 600.0 & 65\% \\ \hline
	\end{tabular}}
%	\begin{tablenotes}
%	\item[1] The test of this algorithm was performed on a $1.8$ times faster 
%machine, a Core i7 920.
%	\end{tablenotes}
\end{threeparttable}
}
\medskip
\caption{Numerical results for 79 realistic instances}
\label{table:Ergebnisse:opt}
\end{table}

\par
In Table \ref{table:Ergebnisse:opt} we summarize the results. More details about 
each instance can be found in Tables \ref{tbl:details:1} and \ref{tbl:details:2}.
We can see that the
relative gaps of all algorithms have the same order of magnitude. The best solutions
are found by the Branch-and-Bound Decomposition, but the improvement of the relative
gap in comparison to the Bid Cut Heuristic is not significant. 
The best suited algorithm in terms of time consumption and quality of the solution
is the Bid Cut Heuristic that finds very good bid selections in about 4 seconds.

\section{Concluding remarks}\label{sec:summary}

Our proposed approach works very well in practice and we were able to
derive desired properties from the formulation as an optimization
problem. A possible route for improvement is to replace the bid
cut by an infeasibility cut of the 
generalized Benders decomposition \cite{geof1972,chu2004}. The main
difficulty here is to find a non-trivial cut that is strong enough to work as fast
as the bid cut. Otherwise the performance of the algorithm for day-ahead 
market clearing might not be sufficient enough. For example the simple
\emph{exact bid cut} (cf. Definition \ref{def:no:good}) used for the Branch-and-Cut
Decomposition is a valid Benders infeasibility cut, but the resulting algorithm is 
too slow.

A nice additional consequence of the presented model are the
transparent pricing rules. In fact, for a given bid selection
optimality can be checked easily: A market participant could
collect all relevant information to check the flow 
price condition. It is clear that a solution will also satisfy the filling 
condition, so by applying Theorem \ref{satz:spotpreiseigenschaften} the participant 
knows that the \surplus{} is maximal for the actual bid selection.
%In this case the market is in an equilibrium. 
This transparency is crucial for market confidence and liquidity. 

We would like to mention that the flow price condition is subject to change in the 
future as it will include \emph{deadbands}: losses during transmission
via an interconnector. While this change will alter the optimality
conditions of our model, our model can be adjusted to incorporate deadbands.

\section*{Acknowledgments}
We would like to thank Herbert Nachbagauer and Adrian Krion for the insightful
discussions, as well as the Deutsche 
B\"orse Systems team for supporting us. We also would 
like to thank European Market Coupling Company for their cooperation and for providing 
us with invaluable data and information.
Furthermore we want to thank the anonymous referees for the valuable comments
and discussions.

%%%%%%%%%%%%%%%%%%%%%%
%%%% Bibtex
%%%%%%%%%%%%%%%%%%%%%%

%\newpage
\urlstyle{same}
\renewcommand\urlstyle[1]{}          % no typewriter in \path commands (e.g., doi-url)
\renewcommand\tt{}                   % no typewriter anywhere else (e.g., arXiv-url)
\bibliographystyle{abbrvnat}
\bibliography{literatur}

%%%%%%%%%%%%%%%%%%%%%%
%%%% Appendix
%%%%%%%%%%%%%%%%%%%%%%

%\newpage
\appendices

\section{Simplifications of the flow price condition}\label{apx:fpc}
\begin{proof}[Proof of eq. \eqref{fpc:simple:1}.]
Suppose that there is no ramping on any interconnector. Then we have 
$\tauramp=\infty$ for all \cons\ so that with \eqref{efp:4} and \eqref{efp:6} 
it follows that $\rhoub\ct=\rholb\ct=0$ for all \cons\ and \hours. Equation 
\eqref{efp:2} before simplifies to
\begin{equation}\label{A1}
	(\muub\ct-\mulb\ct)+\prt-\pst = 0.
\end{equation}
We also know that $\tauct<\tauub$ implies $\muub\ct=0$ and that $\tauct>\taulb$ 
implies $\mulb\ct=0$. Together with \eqref{A1} we therefore conclude
\[
\forall c=\rs\in C, \hours\qquad\pst-\prt=(\muub\ct-\mulb\ct)\begin{cases}
	\leq 0&, \tauct<\tauub\\
	\geq 0&, \taulb<\tauct
\end{cases}.
\]
Using this equation we can interpret the term $\muub\ct-\mulb\ct$ as the
\emph{shadow price of available transmission capacity}: an additional
unit of transmission capacity increases the \surplus\ by
$|\muub\ct-\mulb\ct|$ euro (on a sufficiently small interval).
\end{proof}

\begin{figure}[pb]
	\centering
	\begin{minipage}[c]{12cm}
\includegraphics[width=8cm]{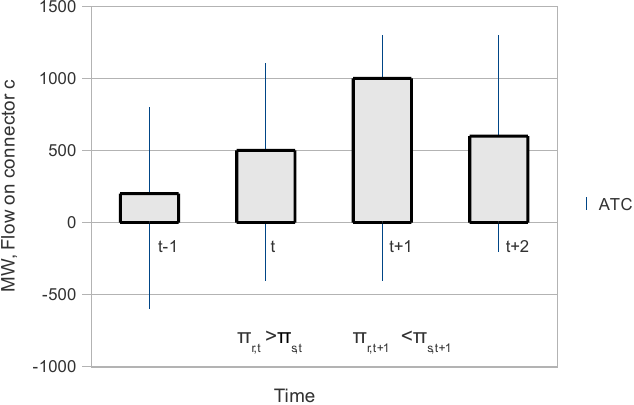}
	\centering
	\end{minipage}
	\caption{Active ramping condition in forward direction, ramp rate = 500 MW}
	\label{fig:active_ramping}
\end{figure}

\begin{proof}[Proof of eq. \eqref{fpc:simple:2}.]
Suppose that there is ramping on the interconnector \cons\ and that the ramping 
condition \nbref{nb:ramprate} is only binding in the two consecutive
time slots $t$ and 
$t+1$. Suppose further that there is still available transfer capacity in that time window.
In other words we have $\taulb<\tauct<\tauub$ and 
$\taulbp<\tau\ctp<\tauubp$. Together with equation \eqref{efp:3} and \eqref{efp:5}
it follows that $\muub\ct=\mulb\ct=\muub\ctp=\mulb\ctp=0$.
We assume that the ramping condition is active 
only between hour $t$ and $t+1$, and that it is active in forward direction (see 
Figure~\ref{fig:active_ramping}). This yields that 
$-\tauramp<(\tau\ctp-\tau\ct)=\tauramp$ and for all $u\in T$ with $u\neq t$ we have 
$-\tauramp<(\tau_{c,u}-\tau_{c,u-1})<\tauramp$. Thus $\rhoub_{c,u}$ and 
$\rholb_{c,u}$ vanish, except for $\rhoub\ctp$.
From equation \eqref{efp:2} and \eqref{efp:1} we get
\begin{align*}
\prt-\pst=&\,\rhoub\ctp\geq 0\;\und\\
\p_{r,t+1}-\p_{s,t+1}=&-\rhoub\ctp\leq 0.
\end{align*}
These statements together provide
$-(\prt-\pst)=-\rhoub\ctp=(\p_{r,t+1}-\p_{s,t+1})$, $\prt>\pst$, and
$\p_{r,t+1}<\p_{s,t+1}$. We can now interpret the variables $\rho$ 
as the shadow prices for ramping, in short the \emph{ramping prices}. Here 
the \surplus\ could by increased by $\rhoub\ctp$ euro if the ramp rate would be 
increased by one unit (again on a sufficiently small interval).
\end{proof}

\begin{proof}[Proof of eq. \eqref{fpc:simple:3a} and \eqref{fpc:simple:3b}.] 
 See also \cite[Theorem 4.4]{muel2009}.
Let the flow-price pair $(\tau,\p)$ satisfy the filling condition, let $c=\rs\in C$, 
\hours, and $\prt<\pst$. Assume that $\tauct<\tauub$ and 
$\tauct-\tau_{c,t-1}<\tauramp$ [, and $\tauct-\tau_{c,t+1}<\tauramp$],
then we have $\muub_{c,t}=0$, $\rhoub_{c,t}=0$ [, and $\rholb_{c,t+1}=0$]. 
Consequently we obtain with \eqref{efp:2} and \eqref{efp:1}
\[
(\underbrace{\muub\ct}_{=0}-\underbrace{\mulb\ct}_{\geq0})
+(\underbrace{\rhoub\ct}_{=0}-\underbrace{\rholb\ct}_{\geq0})
+[\underbrace{\rholb\ctp}_{=0}-\underbrace{\rhoub\ctp}_{\geq0}]
+\underbrace{\prt-\pst}_{<0}=0
\]
which is a contradiction. Implication \eqref{fpc:simple:3b} can be proven 
in a similar fashion.
\end{proof}

\section{MIP formulation of the filling condition}\label{apx:filling:condition}
The filling condition can be easily formulated using binary auxiliary variables. For 
each linear bid curve $h\in H:=\{0,\dots,n\}$ we introduce a binary variable $\gamma_h\in\bin$.
Without loss of generality we assume that the linear bid curves are sorted at first by areas and 
hours, and then by descending prices, that is,
\[
\forall\areas,\hours,\hourly\setminus\max \Hat\qquad 
  m_h(\dhub)\ge m_{h+1}(\underline\delta_{h+1}).
\]
% *****
% I switched from ascending prices to descending prices, such that the sorting is
% consistent with Figure 3.2 and Figure 4.1
Recall that $m_h(\dhlb)\ge m_h(\dhub)$.
The filling condition can now be formulated like this:
\[
\forall\areas,\hours,\hourly\setminus\max \Hat\qquad
% FOR ASCENDING PRICES:
%  \frac{\delta_h-\dhlb}{\dhub-\dhlb}
%    \le\gamma_h\le
%  \frac{\delta_{h+1}-\underline\delta_{h+1}}{\overline\delta_{h+1}-\underline\delta_{h+1}}\;.
% FOR DESCENDING PRICES:
  \frac{\delta_h-\dhlb}{\dhub-\dhlb}
    \ge\gamma_h\ge
  \frac{\delta_{h+1}-\underline\delta_{h+1}}{\overline\delta_{h+1}-\underline\delta_{h+1}}\;.
\]
Note that the presented algorithms do not rely on the inclusion of this condition. In 
fact, including it does not impede the solving process. The slowdown
arises due to the inclusion of 
the filling condition \emph{combined with} the price conditions of
combinatorial bids. This is due to having also to formulate the flow
price condition by using a complicated MIP formulation and this
condition cannot be solved fast enough by standard solvers.

\section{Presolving}
Before starting the optimization process, we can compute upper and lower bounds
for the prices. For every area and hour we can determine the intersection of
the hourly net curve with the price axis. This price would be the clearing 
price if no block bid would be executed and no power would be imported or
exported. Now we assume that we have to import as much additional power as 
possible and that all supply block bids are executed. The additional power
must be bought by hourly bids, thus the price decreases along the hourly net 
curve until all additional power is consumed by hourly bids. This price 
represents a lower bound for the specific area and hour. Upper bounds can be 
determined similarly.

Formula \eqref{eq:price:bounds} is directly derived from the clearing 
constraint \nbref{nb:clearing} and shows the upper and lower quantity 
bounds of the hourly net curves. The price bounds can be derived from these
quantity bounds.
\begin{equation}\label{eq:price:bounds}
\sum_\hourly\delta_h
\begin{cases}
\leq -\sum_\blocka\min\{0, q_{b,t}\}
		-\sum_\flexa\min\{0, q_f\}
		+\sum_\cin\tauub - \sum_\cout\taulb\\
\geq -\sum_\blocka\max\{0, q_{b,t}\}
		-\sum_\flexa\max\{0, q_f\}
		+\sum_\cin\taulb - \sum_\cout\tauub
\end{cases}
\end{equation}
With the help of these bounds the execution state of the blocks and hourly bids
that always realize a profit or incur a loss is given. The bounds can also be 
used to cut off infeasible classes of bid selections in the branch-and-cut
tree. We implemented theses techniques, but unfortunately the solving process 
could not be enhanced significantly.

\begin{table}[p]
\caption{Detailed results, test cases 1-40. Computing times are
measured in milliseconds.}\makebox[\textwidth][c]{\centering\Times{7.5}{
\begin{tabular}{|rrr|rrr|rrr|rrr|}
\hline
\multicolumn{ 1}{|l|}{\textbf{\#}}& \multicolumn{ 1}{c|}{\textbf{Comb.}} & \multicolumn{ 1}{c|}{\textbf{Linear}} & \multicolumn{ 3}{c|}{\textbf{Solver 2.1.2}} & \multicolumn{ 3}{c|}{\textbf{Bid Cut Heuristic}} & \multicolumn{ 3}{c|}{\textbf{BnC Decomposition}} \\ 
\multicolumn{ 1}{|l|}{} & \multicolumn{ 1}{c|}{\textbf{Bids}} & \multicolumn{ 1}{c|}{\textbf{Bid Curves}} & \multicolumn{1}{l}{\textbf{Comp. Time}} & \multicolumn{1}{l}{\textbf{Rel. Gap}} & \multicolumn{1}{l|}{\textbf{PRBs}} & \multicolumn{1}{l}{\textbf{Comp. Time}} & \multicolumn{1}{l}{\textbf{Rel. Gap}} & \multicolumn{1}{l|}{\textbf{PRBs}} & \multicolumn{1}{l}{\textbf{Comp. Time}} & \multicolumn{1}{l}{\textbf{Rel. Gap}} & \multicolumn{1}{l|}{\textbf{PRBs}} \\
\hline\hline
1 & 742 & 34536 & 299000 & 0 & 4 & 2632 & 0 & 4 & 6621 & 0 & 4 \\ 
2 & 727 & 34726 & * 821000 & 4.00E-10 & 8 & 3101 & 0 & 7 & 600565 & 0 & 7 \\ 
3 & 792 & 35227 & 208000 & 3.19E-09 & 5 & 2017 & 3.19E-09 & 5 & 600562 & 3.19E-9 & 5 \\ 
4 & 750 & 34733 & 476000 & 4.77E-08 & 10 & 2621 & 4.61E-08 & 7 & 600334 & 4.61E-8 & 7 \\ 
5 & 685 & 32613 & 190000 & 4.64E-08 & 10 & 1881 & 0 & 6 & 4929 & 0 & 6 \\ 
6 & 630 & 30065 & 350000 & 0 & 3 & 1086 & 0 & 3 & 3992 & 0 & 3 \\ 
7 & 775 & 34516 & * 883000 & 1.34E-08 & 8 & 2615 & 1.31E-08 & 7 & 600251 & 1.31E-8 & 7 \\ 
8 & 772 & 33545 & 243000 & 2.22E-07 & 8 & 2260 & 1.70E-07 & 6 & 600339 & 1.70E-7 & 6 \\ 
9 & 542 & 34876 & 578000 & 2.30E-07 & 5 & 2528 & 0 & 1 & 12653 & 0 & 1 \\ 
10 & 595 & 31226 & 576000 & 1.99E-08 & 6 & 2199 & 0 & 5 & 53924 & 0 & 5 \\ 
11 & 585 & 31888 & 116000 & 2.51E-06 & 17 & 10102 & 9.49E-07 & 5 & 600247 & 8.85E-7 & 5 \\ 
12 & 598 & 33262 & 396000 & 4.51E-07 & 8 & 1968 & 0 & 3 & 9831 & 0 & 3 \\ 
13 & 435 & 29590 & 112000 & 2.74E-06 & 9 & 9970 & 1.14E-06 & 7 & 600213 & 1.14E-6 & 7 \\ 
14 & 558 & 31959 & 460000 & 3.63E-06 & 18 & 4531 & 1.84E-06 & 5 & 600232 & 1.84E-6 & 5 \\ 
15 & 593 & 35384 & 399000 & 0 & 2 & 2143 & 0 & 2 & 7845 & 0 & 2 \\ 
16 & 583 & 36302 & 248000 & 3.29E-06 & 9 & 2868 & 2.83E-06 & 2 & 600339 & 2.83E-6 & 2 \\ 
17 & 585 & 37046 & 602000 & 3.78E-07 & 11 & 5676 & 1.88E-07 & 5 & 600272 & 1.88E-7 & 5 \\ 
18 & 155 & 18195 & 111000 & 0 & 1 & 493 & 0 & 1 & 1114 & 0 & 1 \\ 
19 & 149 & 18287 & 61000 & 0 & 1 & 543 & 0 & 1 & 1135 & 0 & 1 \\ 
20 & 453 & 30471 & 141000 & 2.29E-07 & 9 & 1235 & 2.20E-07 & 2 & 600178 & 2.20E-7 & 2 \\ 
21 & 201 & 30027 & 123000 & 3.21E-06 & 12 & 1871 & 0 & 3 & 6824 & 0 & 3 \\ 
22 & 675 & 30027 & 149000 & 1.73E-06 & 12 & 1824 & 1.50E-06 & 5 & 600225 & 1.50E-6 & 5 \\ 
23 & 701 & 31559 & 135000 & 5.85E-06 & 27 & 34961 & 5.92E-07 & 12 & 600210 & 5.92E-7 & 12 \\ 
24 & 728 & 29902 & 440000 & 3.43E-06 & 4 & 1334 & 3.43E-06 & 4 & 600183 & 3.43E-6 & 4 \\ 
25 & 651 & 30196 & 268000 & 3.24E-06 & 11 & 2366 & 3.11E-06 & 7 & 600216 & 3.11E-6 & 7 \\ 
26 & 642 & 28589 & 162000 & 5.18E-06 & 4 & 2090 & 5.16E-06 & 3 & 600237 & 5.16E-6 & 3 \\ 
27 & 622 & 31692 & 342000 & 3.40E-06 & 21 & 2720 & 3.25E-06 & 7 & 600177 & 3.25E-6 & 7 \\ 
28 & 604 & 31535 & 268000 & 1.08E-07 & 4 & 1551 & 2.59E-08 & 3 & 600204 & 2.59E-8 & 3 \\ 
29 & 707 & 33882 & 611000 & 5.34E-06 & 7 & 3430 & 5.33E-06 & 6 & 600201 & 5.33E-6 & 6 \\ 
30 & 665 & 33821 & 249000 & 1.04E-08 & 2 & 1656 & 1.04E-08 & 2 & 600285 & 1.04E-8 & 2 \\ 
31 & 684 & 33975 & 609000 & 2.39E-05 & 68 & 62814 & 5.04E-06 & 12 & 600231 & 5.04E-6 & 12 \\ 
32 & 599 & 27453 & 181000 & 8.31E-06 & 6 & 2597 & 8.32E-06 & 6 & 600243 & 8.31E-6 & 6 \\ 
33 & 420 & 27247 & 176000 & 9.49E-06 & 3 & 1227 & 9.49E-06 & 3 & 600222 & 9.49E-6 & 3 \\ 
34 & 163 & 24801 & 46000 & 0 & 0 & 802 & 0 & 0 & 2620 & 0 & 0 \\ 
35 & 142 & 26142 & 99000 & 2.40E-08 & 1 & 1020 & 0 & 1 & 3090 & 0 & 1 \\ 
36 & 618 & 31584 & 362000 & 1.57E-07 & 12 & 2189 & 2.23E-07 & 7 & 600226 & 1.23E-7 & 7 \\ 
37 & 249 & 27231 & 94000 & 0 & 1 & 1003 & 0 & 1 & 4403 & 0 & 1 \\ 
38 & 609 & 29589 & 152000 & 3.78E-07 & 27 & 2279 & 4.41E-09 & 11 & 600224 & 4.41E-9 & 11 \\ 
39 & 297 & 30596 & 102000 & 7.50E-07 & 7 & 1342 & 0 & 1 & 5729 & 0 & 1 \\ 
40 & 318 & 30068 & 247000 & 0 & 2 & 1681 & 0 & 2 & 6123 & 0 & 2 \\
\hline
\end{tabular}
}}\label{tbl:details:1}
\medskip\\
*) Computing time was extended in this test so that the
inclusion phase could be finished.
\end{table}

\begin{table}[p]
\caption{Detailed results, test cases 41-79. Computing times are
measured in milliseconds.}\makebox[\textwidth][c]{\centering\Times{7.5}{
\begin{tabular}{|rrr|rrr|rrr|rrr|}
\hline
\multicolumn{ 1}{|l|}{\textbf{\#}}& \multicolumn{ 1}{c|}{\textbf{Comb.}} & \multicolumn{ 1}{c|}{\textbf{Linear}} & \multicolumn{ 3}{c|}{\textbf{Solver 2.1.2}} & \multicolumn{ 3}{c|}{\textbf{Bid Cut Heuristic}} & \multicolumn{ 3}{c|}{\textbf{BnC Decomposition}} \\ 
\multicolumn{ 1}{|l|}{} & \multicolumn{ 1}{c|}{\textbf{Bids}} & \multicolumn{ 1}{c|}{\textbf{Bid Curves}} & \multicolumn{1}{l}{\textbf{Comp. Time}} & \multicolumn{1}{l}{\textbf{Rel. Gap}} & \multicolumn{1}{l|}{\textbf{PRBs}} & \multicolumn{1}{l}{\textbf{Comp. Time}} & \multicolumn{1}{l}{\textbf{Rel. Gap}} & \multicolumn{1}{l|}{\textbf{PRBs}} & \multicolumn{1}{l}{\textbf{Comp. Time}} & \multicolumn{1}{l}{\textbf{Rel. Gap}} & \multicolumn{1}{l|}{\textbf{PRBs}} \\
\hline\hline
41 & 253 & 27288 & 158000 & 6.43E-08 & 5 & 1062 & 0 & 2 & 4626 & 0 & 2 \\ 
42 & 604 & 29969 & 358000 & 7.12E-05 & 2 & 1459 & 7.12E-05 & 2 & 600261 & 7.12E-5 & 2 \\ 
43 & 644 & 32188 & 361000 & 1.96E-07 & 9 & 1658 & 2.89E-08 & 4 & 600246 & 2.89E-8 & 4 \\ 
44 & 683 & 32477 & 232000 & 0 & 2 & 1387 & 0 & 2 & 4125 & 0 & 2 \\ 
45 & 695 & 33131 & 608000 & 1.64E-06 & 7 & 6630 & 1.62E-06 & 4 & 600293 & 1.62E-6 & 4 \\ 
46 & 576 & 28686 & 142000 & 1.07E-06 & 22 & 6151 & 3.83E-07 & 7 & 600215 & 3.83E-7 & 7 \\ 
47 & 476 & 26189 & 181000 & 1.27E-09 & 2 & 752 & 0 & 1 & 1862 & 0 & 1 \\ 
48 & 655 & 29370 & 384000 & 3.90E-08 & 10 & 1666 & 0 & 8 & 16594 & 0 & 8 \\ 
49 & 673 & 32226 & 396000 & 1.37E-05 & 12 & 3937 & 1.34E-05 & 7 & 600198 & 1.34E-5 & 7 \\ 
50 & 667 & 33513 & 349000 & 6.50E-07 & 11 & 1817 & 0 & 1 & 8442 & 0 & 1 \\ 
51 & 738 & 32883 & 621000 & 2.36E-07 & 12 & 2499 & 3.08E-08 & 5 & 600220 & 3.08E-8 & 5 \\ 
52 & 628 & 29959 & 243000 & 4.84E-07 & 26 & 1508 & 2.96E-07 & 3 & 600195 & 2.96E-7 & 3 \\ 
53 & 439 & 29670 & 324000 & 0 & 6 & 1607 & 0 & 6 & 3684 & 0 & 6 \\ 
54 & 705 & 33437 & 604000 & 1.66E-06 & 18 & 2641 & 2.72E-09 & 5 & 600240 & 2.72E-9 & 5 \\ 
55 & 691 & 33609 & 623000 & 3.53E-07 & 10 & 2148 & 0 & 4 & 6236 & 0 & 4 \\ 
56 & 696 & 34061 & 494000 & 6.18E-08 & 6 & 2210 & 2.70E-08 & 4 & 600250 & 2.70E-8 & 4 \\ 
57 & 692 & 33364 & 158000 & 9.78E-07 & 14 & 3282 & 7.97E-07 & 5 & 600246 & 7.97E-7 & 5 \\ 
58 & 644 & 31157 & 254000 & 9.13E-07 & 15 & 3023 & 7.30E-07 & 6 & 600230 & 7.30E-7 & 6 \\ 
59 & 562 & 30123 & 303000 & 1.05E-08 & 5 & 1388 & 0 & 4 & 5158 & 0 & 4 \\ 
60 & 603 & 31833 & 567000 & 2.04E-09 & 7 & 2545 & 2.04E-09 & 7 & 600225 & 2.05E-9 & 7 \\ 
61 & 660 & 34087 & 620000 & 6.46E-07 & 18 & 2259 & 1.89E-07 & 4 & 600293 & 1.89E-7 & 4 \\ 
62 & 682 & 35228 & 611000 & 5.73E-08 & 14 & 2229 & 3.25E-09 & 5 & 600281 & 3.25E-9 & 5 \\ 
63 & 675 & 35128 & 345000 & 0 & 2 & 1561 & 0 & 2 & 5709 & 0 & 2 \\ 
64 & 713 & 34921 & 362000 & 1.62E-07 & 8 & 2180 & 1.30E-07 & 5 & 600326 & 1.30E-7 & 5 \\ 
65 & 627 & 32606 & 203000 & 1.78E-06 & 21 & 1607 & 2.79E-08 & 5 & 600205 & 2.79E-8 & 5 \\ 
66 & 680 & 35368 & 389000 & 0 & 5 & 2090 & 0 & 4 & 7680 & 0 & 4 \\ 
67 & 674 & 34825 & 298000 & 6.11E-06 & 12 & 7555 & 5.75E-06 & 8 & 601153 & 5.75E-6 & 8 \\ 
68 & 711 & 33534 & 610000 & 4.63E-07 & 24 & 1981 & 0 & 5 & 5924 & 0 & 5 \\ 
69 & 730 & 32839 & 626000 & 1.12E-06 & 20 & 5031 & 7.99E-07 & 9 & 600253 & 7.99E-7 & 9 \\ 
70 & 697 & 31728 & 341000 & 6.98E-07 & 23 & 1954 & 4.81E-08 & 4 & 600186 & 4.81E-8 & 4 \\ 
71 & 737 & 33662 & 319000 & 2.36E-07 & 22 & 2816 & 8.64E-09 & 10 & 600230 & 8.64E-9 & 10 \\ 
72 & 708 & 34162 & 497000 & 4.22E-07 & 9 & 1850 & 3.43E-07 & 6 & 600247 & 3.43E-7 & 6 \\ 
73 & 765 & 34315 & 403000 & 7.74E-07 & 7 & 4061 & 4.81E-07 & 7 & 600262 & 5.08E-8 & 5 \\ 
74 & 758 & 36456 & * 745000 & 6.83E-07 & 18 & 3424 & 5.15E-07 & 7 & 600287 & 5.15E-7 & 7 \\ 
75 & 542 & 34437 & 631000 & 1.49E-07 & 4 & 2029 & 0 & 1 & 7967 & 0 & 1 \\ 
76 & 595 & 30812 & 610000 & 3.31E-06 & 18 & 2509 & 0 & 3 & 600227 & 0 & 3 \\ 
77 & 585 & 31593 & 410000 & 1.57E-06 & 13 & 13883 & 1.30E-06 & 6 & 600248 & 1.30E-6 & 6 \\ 
78 & 598 & 32966 & 479000 & 4.80E-07 & 12 & 2093 & 0 & 3 & 51173 & 0 & 3 \\ 
79 & 558 & 31525 & 605000 & 3.87E-06 & 16 & 21331 & 1.52E-06 & 6 & 600212 & 1.52E-6 & 6 \\\hline\hline
\textit{avg.} & \textit{595} & \textit{31692} & \textit{361380} & \textit{2.58E-06} & \textit{10.76} & \textit{4114} & \textit{1.93E-06} & \textit{4.65} & \textit{390808} & \textit{1.92E-6} & \textit{4.62} \\\hline
\end{tabular}
}}\label{tbl:details:2}
\medskip\\
*) Computing time was extended in this test so that the
inclusion phase could be finished.
\end{table}

\end{document}

%% file: texDefinitions.tex
\global\long\def\poly{\mathrm{poly}}

%% file: commands.tex
\newcommand{\surplus}{economic surplus}

\newcommand{\epex}{EPEX Spot}

\newcommand{\mycommand}[2]{\ifdefined#1\renewcommand{#1}{#2}\else\newcommand{#1}[0]{#2}\fi}

\newcommand{\comment}[1]{}

\newcommand{\R}{\mathbb R}
\newcommand{\Q}{{\mathbb R}}
\newcommand{\bin}{\{0,1\}}
\newcommand{\eur}{\geneuronarrow}%\officialeuro}

\mycommand{\mit}{\text{ with }}
\newcommand{\und}{\text{ and }}

\newcommand{\fuer}{\text{ for }}

\newcommand{\dd}{\mathrm{d}}

\newcommand{\pricerange}{\R}
\renewcommand{\Hat}{{H_{a,t}}}

\newcommand{\hourly}{{h\in H_{a,t}}}

\newcommand{\blocks}{{\ensuremath{b\in B}}}
\newcommand{\blocka}{{\ensuremath{b\in B_a}}}
\newcommand{\bt}{_{b,t}}
\newcommand{\flex}{{\ensuremath{f\in F}}}
\newcommand{\flexa}{{\ensuremath{f\in F_a}}}
\newcommand{\hours}{{\ensuremath{t\in T}}}
\newcommand{\areas}{{\ensuremath{a\in A}}}
\newcommand{\rs}{(r,s)}
\newcommand{\cons}{{\ensuremath{{c\in C}}}}
\newcommand{\cout}{{c\in C^+_a}}
\newcommand{\cin}{{c\in C^-_a}}
\newcommand{\tauct}{\tau_{c,t}}
\newcommand{\tauub}{\overline\tau_{c,t}}
\newcommand{\taulb}{\underline\tau_{c,t}}
\newcommand{\tauubp}{\overline\tau_{c,t+1}}
\newcommand{\taulbp}{\underline\tau_{c,t+1}}
\newcommand{\tauramp}{\widetilde{\tau}_{c}}

\newcommand{\ct}{_{c,t}}
\newcommand{\ctp}{_{c,t+1}}
\newcommand{\ctm}{_{c,t-1}}

\newcommand{\dph}{\Delta p_h}

\newcommand{\pft}{\varphi_{f,t}}

\newcommand{\at}{_{a,t}}

\newcommand{\p}{\ensuremath{\pi}}
\newcommand{\pat}{\ensuremath{\pi_{a,t}}}
\newcommand{\prt}{\pi_{r,t}}
\newcommand{\pst}{\pi_{s,t}}
\newcommand{\pmin}{\underline p}
\newcommand{\pmax}{\overline p}
\mycommand{\poly}{\ensuremath{\mathcal{P}}}
\newcommand{\valP}{\ensuremath{(\p,\beta,\varphi,\delta,\tau)}}

\newcommand{\muub}{\overline\mu}
\newcommand{\mulb}{\underline\mu}
\newcommand{\rhoub}{\overline\rho}%{\widetilde\rho}
\newcommand{\rholb}{\underline\rho}%{\undertilde\rho}

\newcommand{\valp}{\p^{*}}
\newcommand{\valpat}{\valp_{a,t}}

\newcommand{\valblockgain}{\sum_\hours (p_b-\valpat)q_{b,t}}

\newcommand{\valflexgain}{(p_f-\valpat)\sgn(q_f)}

\newcommand{\nbref}[1]{(\ref{#1})}
\newcommand{\qpflow}{\eqref{QPFlow}}
\newcommand{\qprelax}{\eqref{QPRelax}}
\newcommand{\fixflow}{\eqref{FixFlow}}

\newcommand{\qpbidcut}{\eqref{QPBidCut}}
\newcommand{\qpprice}{\eqref{QPPrice}}
\ifdefined\Satz\renewcommand{\Satz}[1]{Satz \ref{satz:#1}}\else\newcommand{\Satz}[1]{Satz \ref{satz:#1}}\fi
\newcommand{\tp}{^{\top}}

\newcommand{\valbidcut}{\ensuremath{(\valbeta,\valphi,\valdelta,\valtau)}}
\newcommand{\valqpflow}{\ensuremath{(\valp,\valbeta,\valphi,\valdelta,\valtau)}}
\newcommand{\varqpflow}{\ensuremath{(\p,\beta,\varphi,\delta,\tau)}}

\newcommand{\cplex}{{IBM CPLEX}}

\newcommand{\phvonq}{m_h}%\mathbf{}
\newcommand{\phvond}{\phvonq(\delta_h)}

\newcommand{\dhub}{\overline\delta_h}
\newcommand{\dhlb}{\underline\delta_h}

\newcommand{\phd}{p_h-\dph\delta_h}

\newcommand{\deltawithinbounds}{$\dhlb\le\delta_h\le\dhub$ for $h\in H$}

\newcommand{\valtau}{\tau^{*}}
\newcommand{\valtauct}{\tau_{c,t}^*}
\newcommand{\valdelta}{\delta^*}
\newcommand{\valbeta}{\beta^*}
\newcommand{\valphi}{\varphi^*}

\newcommand{\valrelax}{\ensuremath{(\valdelta,\valtau)}}

\newcommand{\vub}{\overline{v}}
\newcommand{\vlb}{\underline{v}}
\newcommand{\gebotswahl}{\ensuremath{(\valbeta,\valphi)}}

\newcommand{\Bloss}{B_\text{loss}}
\newcommand{\Floss}{F_\text{loss}}
\newcommand{\Blossi}{B_\text{loss}^i}
\newcommand{\Flossi}{F_\text{loss}^i}
\newcommand{\cut}{{\rm Cut}}